\documentclass[12pt,a4paper,twoside]{article}
\usepackage[T1]{fontenc}
\usepackage[cp1250]{inputenc}

\usepackage[polish,english]{babel} 
\usepackage{graphicx}
\usepackage[intlimits]{amsmath}
\usepackage{amsthm}
\usepackage{amsfonts}

\usepackage{stmaryrd}
         \newcommand{\dow}{{\noindent\it Proof.~}}
\newcommand{\calka}{\int_{\Omega}}
\newcommand{\calkam}{\int_{\Omega^{-}}}
\newcommand{\calkap}{\int_{\Omega^{+}}}
\newcommand{\calkad}{\int_{\partial\Omega}}
\newcommand{\calkadm}{\int_{\partial\Omega^{-}}}
\newcommand{\calkadp}{\int_{\partial\Omega^{+}}}
\newcommand{\calkat}{\int_{0}^{T}}
\newcommand{\dt}{\frac{1}{\triangle t}}
\newcommand{\Dt}{\partial_{t}}
\newcommand{\dyv}{\mathrm{\ div}}
\newcommand{\rot}{\mathrm{rot}}

\newtheorem{lemat}{Lemma}
\newtheorem{prop}[lemat]{Proposition}
\newtheorem{twr}[lemat]{Theorem}
\newtheorem{definition}{Definition}
\newtheorem{remark}[lemat]{Remark}

\title{Analysis of  semidiscretization of the compressible Navier-Stokes equations.}
\author{Ewelina Kamińska}

\begin{document}
\maketitle
\normalsize
\begin{center}
{\small Institute of Applied Mathematics and Mechanics\\
University of Warsaw\\ul.Banacha 2, 02-097 Warszawa, Poland\\
{\sc E-mail: e.kaminska@mimuw.edu.pl}}
\end{center}
\noindent{\bf Abstract:} The objective of this work is to present the existence result for the evolutionary  compressible Navier-Stokes equations via time discretization. We consider the two-dimensional case with  slip boundary conditions. First, the existence of  weak solution for a fixed  time step $\triangle t>0$ is presented and then the limit passage as $\triangle t\rightarrow 0^{+}$ is carried out. The proof is based on a new technique established for the stedy Navier-Stokes equations by Mucha P. and Pokorn\'y M. 2006 {\it Nonlinearity} {\bf19} 1747-1768.

\vspace{0.2cm} \noindent{\it Mathematics Subject
Classification:} 76N10; \vspace{0.2cm} 35Q30

\section{Introduction}
We investigate a system being time discretization of two dimensional Navier-Stokes equations in the isentropic regime
\begin{equation}\label{eq:system}
\begin{array}{r}
\dt\left(\varrho^{k}-\varrho^{k-1}\right)+\dyv(\varrho^{k}v^{k})=0\\
\dt\left(\varrho^{k}v^{k}-\varrho^{k-1}v^{k-1}\right)+\dyv(\varrho^{k}v^{k}\varotimes v^{k})-\mu\Delta v^{k}-(\mu+\nu)\nabla\dyv v^{k}+\nabla\pi(\varrho^{k})=0,
\end{array}
\end{equation}
where $\Omega\subset\mathbb{R}^{2}$ is a fixed domain, $v^{k}:\Omega\rightarrow\mathbb{R}^{2}$- the velocity field, $\varrho^{k}:\Omega\rightarrow\mathbb{R}^{+}_{0}$- the density, $\pi:\mathbb{R}^{+}_{0}\rightarrow \mathbb{R}$- the internal pressure given by the constitutive equation
\begin{equation*}
\pi(\varrho^{k})=(\varrho^{k})^{\gamma},\ \gamma>2.
\end{equation*}
We assume that the walls of $\Omega$ are rigid and that the fluid slips at the boundary
\begin{equation}\label{slip}
\begin{array}{rcl}
v^{k}\cdotp n=0,\quad at\ \partial\Omega\\
n\cdotp T(v^{k},\pi)\cdotp\tau+fv^{k}\cdotp\tau=0\quad at\ \partial\Omega,
\end{array}
\end{equation}
where $T(v^{k},\pi)=2\mu D(v^{k})+(\nu\dyv v-\pi)I$.\\
The conditions (\ref{slip}) are known as the Navier or friction relations which means, that unlike in the case of complete slip of the fluid against the boundary, the friction effects, described by $f\geq0$, may also be present. The customary zero Dirichlet condition may be understood as a special case of the above, when $f\rightarrow\infty$.\\
We will always assume that our initial conditions $\varrho^{0}, v^{0}$ satisfy
\begin{eqnarray*}
\varrho^{0}\geq 0\ a.e.\ in\ \Omega,\quad \varrho^{0}\in L_{\gamma}(\Omega),\\
\varrho^{0}v^{0}\in L_{2\gamma/(\gamma+1)}(\Omega),\quad \varrho^{0}(v^{0})^{2}\in L_{1}(\Omega).
\end{eqnarray*}
The first main goal of this paper is to show that for $\triangle t=const.$ the solutions of such a system exist in a sense of the following definition.
\begin{definition}
We say, the pair of functions $(\varrho^{k},v^{k})\in L_{\gamma}(\Omega)\times W_{2}^{1}(\Omega)$, $v^{k}\cdotp n=0$ at $\partial\Omega$ is a weak solution to (\ref{eq:system})-(\ref{slip}) provided
\begin{equation*}
\calka\varrho^{k}v^{k}\cdot\nabla\varphi\ dx=\dt\calka(\varrho^{k}-\varrho^{k-1})\varphi\ dx,\quad\forall\varphi\in C^{\infty}(\overline{\Omega}),
\end{equation*}
and
\begin{multline*}
\dt\calka(\varrho^{k}v^{k}-\varrho^{k-1}v^{k-1})\varphi\ dx-\calka\varrho^{k}v^{k}\varotimes v^{k}:\nabla\varphi\ dx+2\mu\calka\mathbf{D}(v^{k}):\mathbf{D}(\varphi)\ dx\\+\nu\calka\dyv{v^{k}}\dyv{\varphi}\ dx-\calka\pi(\varrho^{k})\dyv{\varphi}\ dx+\calkad f(v^{k}\cdot\tau)(\varphi\cdot\tau)\ dS=0,\quad\\ \forall\varphi\in C^{\infty}(\overline{\Omega});\ \varphi\cdot n=0\ at\ \partial\Omega.
\end{multline*}
\end{definition}
\noindent The first main result reads as follows.
\begin{twr}\label{theorem1}
Let $\Omega\in C^{2}$ be a bounded domain, $\triangle t=const.$, $\mu>0$, $2\mu+3\nu>0$, $\gamma>2$, $f\geq 0$, $\varrho^{k-1}\geq 0$. Then there exists a weak solution to (\ref{eq:system})-(\ref{slip}) such that
\begin{eqnarray*}
&{}&\varrho^{k}\in L_{\infty}(\Omega)\quad and\ \varrho^{k}\geq 0,\\
&{}&v^{k}\in W^{1}_{p}(\Omega)\quad  \forall p<\infty,\\
&{}&\calka\varrho^{k}dx=\calka\varrho^{k-1}dx,
\end{eqnarray*}
moreover $\|\varrho^{k}\|_{\infty}\leq (\triangle t)^{3/(1-\gamma)}$.
\end{twr}
The result we present here was already stated without requiring any assumption on the smallness of initial data  f.i. in the monograph of Lions  \cite{PLL} for the zero Dirichlet condition when $\Omega$ is bounded and for the whole space. It was used there as a tool in analysing the steady and non-steady cases. The approach presented there or in other works based on the Feireisl idea \cite{NS}, \cite{EF1}, \cite{FN} benefit from the properties of the effective viscous flux. Our technique allows for essential reduction of the number of technical tricks and enables to get the required  $L_{\infty}$ regularity for the density directly at the level of approximate system.\\
The second result refers to a passage to the limit  with length of time interval $\triangle t\rightarrow 0$. We will show that for such a case our solution tends to the weak solution of non-steady compressible Navier-Stokes system with a slip boundary condition:
\begin{equation}\label{eq:system2}
\begin{array}{rlr}
\varrho_{t}+\dyv(\varrho v)=0\quad&in &\Omega\\
(\varrho v)_{t}+\dyv(\varrho v\varotimes v)-\mu\Delta v-(\mu+\nu)\nabla\dyv v+\nabla\pi(\varrho)=0\quad&in&\Omega\\
v\cdotp n=0\quad&at&\partial\Omega\\
n\cdotp T(v,\pi)\cdotp\tau+fv\cdotp\tau=0\quad&at&\partial\Omega,
\end{array}
\end{equation}
in sense of the following definition.
\begin{definition}
We say, the pair of functions $(\varrho,v)\in L_{\infty}(L_{\gamma})\times L_{2}(W^{1}_{2})$, $v\cdot n=0$ at $\partial\Omega$ is a weak solution to (\ref{eq:system2}) provided
\begin{equation*}
\calkat\calka\left(\varrho\varphi_{t}+\varrho v\cdot\nabla\varphi\right)\ dxdt=0, \quad \forall\varphi\in C^{\infty}_{c}([0,T)\times\overline{\Omega}),
\end{equation*}
and
\begin{multline}
\calkat\calka\left(\varrho v\varphi_{t}+\varrho v\varotimes v:\nabla_{x}\varphi+\pi(\varrho)\dyv_{x}{\varphi}\right)\ dxdt=\\
=\calkat\calka\left(2\mu\mathbf{D_{x}}(v):\mathbf{D_{x}}(\varphi)+\nu\dyv_{x}{v}\ \dyv_{x}{\varphi}\right)\ dxdt+\calkat\calkad f(v\cdot\tau)(\varphi\cdot\tau)\ dSdt,\\ \forall\varphi\in C^{\infty}_{c}([0,T)\times\overline{\Omega});\ \varphi\cdot n=0\ at\ \partial\Omega.
\end{multline}
\end{definition}
\noindent  The existence of solutions to the non-steady system is provided by our second main result.
\begin{twr}\label{theorem2}
Under the hypotheses of Theorem \ref{theorem1}, the solution $(\varrho^{k},v^{k})$ converges to $(\varrho,v)$ as $\triangle t\rightarrow 0^{+}$ weakly (weakly${}^{*}$) in $ L_{\infty}(L_{\gamma})\times L_{2}(W^{1}_{2})$. \\
Moreover $\varrho$ belongs to $L_{\gamma+1}(\Omega\times(0,T))$ and the following energy inequality is satisfied for almost all $t\in[0,T]$
\begin{multline*}
\calka\varrho v^{2}(T)dx+\frac{1}{\gamma-1}\calka\varrho^{\gamma}(T)dx+\calkat \calka\left( 2\mu |D(v)|^{2}+\nu(\dyv v)^{2}\right)\ dxdt\\
+\calkat\calkad f(v\cdot \tau)^{2}dx\ dt\leq C(\varrho^{0},v^{0}).
\end{multline*}
\end{twr}
In the following section we will show the existence and uniqueness  of regular solution to the problem being the new $\epsilon-$approximation  scheme for the time-discretized Navier-Stokes equations.  Although the proof is based on the standard fixed-point method, we will precisely present most of steps in view of the fact that our approximation affects the nonlinear term too. Our solution $(\varrho^{k},\ v^{k})$ will be obtained as a weak limit as $\epsilon\rightarrow 0^{+}$ of the sequences $(\varrho^{k}_{\epsilon},\ v^{k}_{\epsilon})$. This limit process will be carried out in Section 3 by using some uniform estimates and the following property of the density sequence 
\begin{equation*}
\lim_{\epsilon\rightarrow 0^{+}}|\{x\in\Omega:\varrho^{k}_{\epsilon}(x)>m\}|=0
\end{equation*}
for $m$ sufficiently large, which enables to show the convergence of the pressure.

\section{Approximation}

In this section we  present a scheme of approximation being a modification of the one introduced by Mucha  Pokorny \cite{MP} for the steady case. It is needed  to investigate  the issue of existence of solutions in the case when the time step ($\triangle t$) is cosnstant and while disposing a sufficient information for the density and velocity at the $k-1$ moment of time. Although for further purposes there is a necessity to keep trace of the dependence on these quantities in almost all estimates. \\
Let
	\begin{equation}\label{notacja1}
	\begin{array}{c}
	\alpha=\dt,\\
	h=\varrho^{k-1},\quad
	\varrho=\varrho^{k},\quad
	v=v^{k},\quad
	g=v^{k-1}.
	\end{array}
	\end{equation}
The objective of this part of work will be then to examine the following approximative system:
\begin{equation}\label{approx}
\begin{array}{r}
\alpha\left(\varrho-hK(\varrho)\right)+\dyv(K(\varrho)\varrho v)-\epsilon\Delta\varrho=0\\
\alpha\left(\varrho v-hg\right) +\dyv(K(\varrho)\varrho v\varotimes v)-\mu\Delta v-(\mu+\nu)\nabla\dyv v+\nabla P(\varrho)+\epsilon\nabla\varrho\nabla v=0\\
\frac{\partial\varrho}{\partial n}=0\quad at\quad \partial\Omega,\\
v\cdot n=0\quad at\quad \partial\Omega,\\
n\cdotp T(v,P(\varrho))\cdotp\tau+fv\cdotp\tau=0\quad at\quad \partial\Omega,\\
\end{array}
\end{equation}
we will write simply $\varrho, v$ istead of $\varrho_{\epsilon}, v_{\epsilon}$ when no confusion can arise. The other denotations are the following:
\begin{equation*}
P(\varrho)=\gamma\int_{0}^{\varrho}s^{\gamma-1}K(s)ds,
\end{equation*}
where
\begin{equation*}
K(\varrho)=\left\{
\begin{array}{ll}
1 &\quad \varrho\leq m_{1},\\
0 &\quad \varrho\geq m_{2},\\
\in(0,1) &\quad \varrho\in(m_{1},m_{2}),
\end{array}
\right.
\end{equation*}
and
\begin{equation*}
K(\cdot)\in C^{1}(\mathbb{R})\quad K'(\varrho)<0\  in\ (m_{1},m_{2}), 
\end{equation*}
for some constants $m_{1},\ m_{2}$. To avoid the difficulties conected with the case when $m_{1}\rightarrow m_{2}$ we set the difference $m_{2}-m_{1}$ to be constant, equal 1.\\
The existence of a regular solution is guaranteed by the theorem.
\begin{twr}\label{twr:approx}
Let $\Omega\in C^{2},\ \epsilon , \varrho_{0}, \dt>0$. Then there exist a regular solution $(\varrho, v)$ to (\ref{approx}), $\varrho\in W^{2}_{p}(\Omega),\ v\in W^{2}_{p}(\Omega)$ for all $p<\infty$.\\
Moreover
\begin{eqnarray}
0\leq\varrho\leq m_{2}\quad in\ \Omega,\label{nonnegativ}\\
\calka \varrho dx\leq\calka h dx.\label{const}
\end{eqnarray}
\end{twr}
\dow
We assume, that $\varrho,v$ are regular solutions to (\ref{approx}) and prove some estimates first, after we go on with the existence.\\
\noindent Step 1. {\it Proof of (\ref{const})}.\\
Integrating  the first equation of (\ref{approx}) over $\Omega$ one gets
 \begin{equation*}
 \alpha\calka(\varrho-hK(\varrho))dx+\calkad K(\varrho)\varrho v\cdot n dS-\epsilon\calkad\frac{\partial\varrho}{\partial n}dS=0,
 \end{equation*}
 the boudary integrals vanish and due to the definition of $K(\cdot)$ we truly have
 \begin{equation*}
 \calka\varrho dx= \calka K(\varrho) h dx\leq\calka hdx.
 \end{equation*}
 \noindent Step 2. {\it Non-negativity of $\varrho$.}  \\
 Assume, that we have $h\geq0$ in $\Omega$, the proof follows by the induction. We integrate first equation of (\ref{approx}) over $\Omega^{-}=\{x\in\Omega : \varrho(x)<0\}$
 	 \begin{equation*}
	 \alpha\calkam(\varrho-K(\varrho)h)dx+\calkadm K(\varrho)\varrho v\cdot n dS-\epsilon\calkadm\frac{\partial\varrho}{\partial n}dS=0,
 	 \end{equation*}
 the first boundary integral vanishes since either $\varrho$ or $v\cdot n$ equals $0$ at  $\partial\Omega^{-}$. Moreover, we know that  $\frac{\partial\varrho}{\partial n}\geq 0$ at $\partial \Omega^{-}$, hence
	 \begin{equation*}
 	 \calkam\varrho dx\geq\calkam K(\varrho)hdx\geq 0,
	 \end{equation*}
 but this leads to conclusion that $|\Omega^{-}|=0$ and consequently $\varrho\geq 0$ in $\Omega$.\\
Step 3. {\it Upper bound for $\varrho$.}\\
Assume that $h\leq m_{2}$. This time we integrate the approximate continuity equation over $\Omega^{+}=\{x\in\Omega : \varrho(x)\geq m_{2}\}$
\begin{equation*}
	 \alpha\calkap(\varrho-K(\varrho)h)dx+\calkadp K(\varrho)\varrho v\cdot n dS-\epsilon\calkadp\frac{\partial\varrho}{\partial n}dS=0,
\end{equation*}
At $\partial\Omega^{+}$ we have $\frac{\partial\varrho}{\partial n}\leq 0$ and either $K(\varrho)$ or $v\cdot n$ equals $0$. Thus, in the similar way as previously, the observation
\begin{equation*}
\calkap\varrho dx\leq m_{2}\calkap K(\varrho)dx\leq 0
\end{equation*}
implies that $\varrho\leq m_{2}$ in $\Omega$.\\
Step 4. {\it Existence.}\\
In accordance with our denotations the proof of existence of approximate solutions is almost identical to the one presented in \cite{MP}.
In the first step we define for $p\in[1,\infty]$:
	\begin{equation*}
	M_{p}=\left\{w\in W^{1}_{p}(\Omega);w\cdot n=0\  at\ \partial\Omega \right\}.
	\end{equation*}
and we claim that the following proposition, which is the analogue of Proposition 3.1. from \cite{MP} holds true.
\begin{prop}\label{S}
Let assumptions of theorem \ref{twr:approx} be satisfied. Then the operator $S:M_{\infty}\rightarrow W^{2}_{p}(\Omega)$, where
\begin{eqnarray*}
S(v)=\varrho,\\
\alpha\varrho +\dyv(K(\varrho)\varrho v)-\epsilon\Delta\varrho=\alpha hK(\varrho)\quad in\quad\Omega\\
\frac{\partial\varrho}{\partial n}=0\quad at\quad\partial\Omega
\end{eqnarray*}
is well defined for any $p<\infty$. Moreover\\
\begin{itemize}
\item $\varrho=S(v)$ satisfy 
\begin{equation*}
\calka\varrho dx\leq \calka hdx.
\end{equation*}
\item If $h\geq 0$ then $\varrho\geq0$ a.e. in $\Omega$.
\item If $\|v\|_{1,\infty}\leq L,\ L>0$ then
\begin{equation}\label{hireg}
\|\varrho\|_{2,p}\leq C(\epsilon,p,\Omega)(1+L)\|h\|_{p},\quad 1<p<\infty.
\end{equation}
\end{itemize}
\end{prop}

 The only difference in the formulation and the proof relates to the fact that $h$ is not a constant parameter any more and that  there appears $g$ instead of $v$. But the assumtion that  the regular solution in $k-1$ moment of time  exist allows to replace the modulus by  the $L_{p}$ norm of $h$.\\
\noindent In the next step of proof of Theorem \ref{twr:approx} we will consider the Lame operator 
	\begin{equation*}
	\mathcal{T}:\ M_{\infty}\rightarrow M_{\infty} 
	\end{equation*}
defined as follows: $w=\mathcal{T}(v)$ is a solution to the problem
	\begin{equation}\label{Lame}
	\begin{array}{r}
	-\mu \Delta w-(\mu+\nu)\nabla\dyv w=\alpha hg-\alpha\varrho v-\dyv(K(\varrho)\varrho v\varotimes v)-\nabla P(\varrho)-\epsilon\nabla\varrho \nabla v=\\
	=F(\varrho,v,h,g)\\
	w\cdot n=0\quad at\quad \partial\Omega,\\
	n\cdotp (2\mu D(w)+\nu\dyv w I)\cdotp\tau+fv\cdotp\tau=0\quad at\quad\ \partial\Omega\\
	\end{array}
	\end{equation}
Employing the Larey-Schauder fixed point theorem for the operator $\mathcal{T}$ we can almost rewrite the proof of analogous fact \cite{NS} or \cite{MP}. The only part that that deserves more careful study is the energy estimate. This in turn together with some information about the pressure $P(\varrho)$ will enable to pass to the limit with the length of time interval $\triangle t$.
First observe that  $(\ref{Lame})_{1}$ with $w=v$ and $\varrho=S(v)$ holds with  a solution itself as a test function, therefore
\begin{multline*}
	\alpha\calka\varrho v^{2}+\calka\dyv(K(\varrho)\varrho v\varotimes v)v-\mu\calka(\Delta v)v-(\mu+\nu)\calka(\nabla\dyv v)v+\calka\nabla P(\varrho)v\\
	+\epsilon\calka\nabla\varrho \nabla v v
	=\alpha \calka hgv.
	\end{multline*}
Next, integrating by parts and using condition on the boundary
	\begin{multline*}
	\alpha\calka\varrho v^{2}+\frac{1}{2}\calka\dyv(K(\varrho)\varrho v)v^{2}+2\mu\calka|D(v)|^{2}+\nu\calka\dyv^{2}v+\calkad f(v\cdotp\tau)^{2}\\
	-\frac{\gamma}{\gamma-1}\calka\dyv(K(\varrho)\varrho v)\varrho^{\gamma-1}-\frac{\epsilon}{2}\calka\Delta\varrho v^{2}
	=\alpha \calka hgv,
	\end{multline*}
including the information contained in $(\ref{approx})_{1}$ one gets
	\begin{multline*}
	\frac{1}{2}\alpha \calka(\varrho+K(\varrho)h)v^2+2\mu\calka|D(v)|^{2}+\nu\calka\dyv^{2}v+\calkad f(v\cdotp\tau)^{2}
	\\+\frac{\gamma}{\gamma-1}\alpha\calka\varrho^{\gamma}-\frac{\gamma}{\gamma-1}\alpha\calka\varrho^{\gamma-1}K(\varrho)h+\gamma\epsilon \calka|\nabla\varrho|^{2}\varrho^{\gamma-2}
	=\alpha \calka hgv,
	\end{multline*}
now we add end substract $\frac{1}{2}\alpha \calka hg^{2}$
	\begin{multline}
	\frac{1}{2}\alpha\calka(\varrho v^2-hg^2)+\frac{1}{2}\alpha\calka h|v-g|^{2}+
	2\mu\calka|D(v)|^{2}+\nu\calka\dyv^{2}v\\+
	\calkad f(v\cdotp\tau)^{2}
	+\frac{\gamma}{\gamma-1}\alpha\calka\varrho^{\gamma}-\frac{\gamma}{\gamma-1}\alpha\calka\varrho^{\gamma-1}K(\varrho)h+
	\gamma\epsilon \calka|\nabla\varrho^{\frac{\gamma}{2}}|^{2}\leq0,
	\end{multline}
next we add and substract $\frac{1}{\gamma-1}\alpha \calka h^{\gamma}$
	\begin{multline}\label{koszmar:approx}
	\frac{1}{2}\alpha\calka(\varrho v^2-hg^2)+\frac{1}{2}\alpha\calka h|v-g|^{2}+
	2\mu\calka|D(v)|^{2}+\nu\calka\dyv^{2}v+
	\calkad f(v\cdotp\tau)^{2}\\
	+\frac{1}{\gamma-1}\alpha\calka\left(\varrho^{\gamma}-h^{\gamma}\right)+\frac{1}{\gamma-1}\alpha\calka\left( (\gamma-1)\varrho^{\gamma}+h^{\gamma}-\gamma\varrho^		{\gamma-1}K(\varrho)h\right)+\gamma\epsilon \calka|\nabla\varrho^{\frac{\gamma}{2}}|^{2}\leq0.
	\end{multline}
Note, that since $\varrho,\ h\geq 0$ and $K(\varrho)\leq 1$ we have that $(\gamma-1)\varrho^{\gamma}+h^{\gamma}-\gamma\varrho^{\gamma-1}K(\varrho)h\geq 0$. \\
Referring to our oryginal denotation we may now sum (\ref{koszmar:approx}) from $k=1$ to $k=n$ and obtain the following bounds:
	\begin{equation*}
	\sup_{0\leq n\leq M}\frac{1}{\gamma-1}\alpha \|\varrho^{n}\|^{\gamma}_{\gamma}+\frac{1}{2}\alpha \|\varrho^{n}(v^{n})^2\|_{1}\leq \frac{1}{\gamma-1}\alpha \|\varrho^{0}\|^		{\gamma}_{\gamma}+\frac{1}{2}\alpha \|\varrho^{0}(v^{0})^2\|_{1} ,
	\end{equation*}
thus
	\begin{equation}\label{osz1}
	\sup_{0\leq n\leq M}\|\varrho^{n}\|^{\gamma}_{\gamma}+\|\varrho^{n}(v^{n})^2\|_{1}\leq C(\varrho^{0},v^{0},\gamma,\Omega),
	\end{equation}
in particular C is independent of $k,\epsilon$ and $\alpha$, moreover
	\begin{equation}\label{osz2}
	\sum_{k=1}^{M}\calka\left[\varrho^{k-1}|v^{k}-v^{k-1}|^{2}+ (\gamma-1)(\varrho^{k})^{\gamma}+(\varrho^{k-1})^{\gamma}-\gamma(\varrho^{k})^{\gamma-1}K(\varrho^{k})\varrho^{k-1}\right]		\leq C
	\end{equation}
with the same constant $C$. The information contained here turns aut to be one of the crutial importance at the second stage of this work while showing that the passage with $\triangle t\rightarrow 0$ gives the solution to the evolutionary case. Namely, since for $\gamma>2$ there exists a positive constant $\delta$, such that
\begin{equation*}
(\gamma-1)(\varrho^{k})^{\gamma}+(\varrho^{k-1})^{\gamma}-\gamma(\varrho^{k})^{\gamma-1}K(\varrho^{k})\varrho^{k-1}\geq\delta|\varrho^{k}-\varrho^{k-1}|^{\gamma},
\end{equation*}
hence (\ref{osz2}) ensures
\begin{equation}\label{osz3}
\sum_{k=1}^{M}|\varrho^{k}-\varrho^{k-1}|^{\gamma}\leq C.
\end{equation}

Additionaly we have
	\begin{equation*}
	\sum_{k=1}^{M}\|Dv^{k}\|^{2}_{2}\leq \alpha  C
	\end{equation*}
and by Korn's inequality 
	\begin{equation}\label{vH1}
	\sum_{k=1}^{M}\|v^{k}\|^{2}_{1,2}\leq \alpha  C
	\end{equation}
here the constant $C$ depends also on $\mu$ and $\nu$.\\
Finally we also get
	\begin{equation}\label{rhoH1}
	\sum_{k=1}^{M}\|\nabla (\varrho^{k})^{\frac{\gamma}{2}}\|^{2}_{2}\leq \frac{\alpha}{\epsilon}C.
	\end{equation}	
This information allows us to repeat the procedure described in \cite{NS}, which together with the Proposition \ref{S} yield the existence of regular solutions, and hence the proof of Theorem \ref{twr:approx} is complete.\\
Apart from the information resulting from the first {\it a priori} estimate , the limit passage requires also some estimates independent  $\epsilon,\ \alpha$ and $m_{2}$. \\
First of them is the estimate for the norm of gradient of the density. Observe that multiplying $(\ref{approx})_{1}$ by $\varrho$ and integrating over $\Omega$ one get
\begin{multline*}
\epsilon\calka|\nabla\varrho|^{2}=\alpha\calka hK(\varrho)\varrho-\alpha\calka \varrho^{2}-\calka K(\varrho)\varrho v\cdot\nabla\varrho\\
\leq \alpha C m_{2}+\calka v\cdot\nabla\left(\int_{0}^{\varrho}K(t) t\  dt \right)=\alpha C m_{2}-\calka \dyv v\left(\int_{0}^{\varrho}K(t) t\  dt \right)\\
\leq \alpha C m_{2}+\calka |\dyv v|\varrho^{2}\leq \alpha C m_{2}+\sqrt{\alpha}C m_{2}^{2}.
\end{multline*}
This means that $\|\nabla\varrho\|_{2}$ may blow up as $\epsilon\rightarrow 0^{+}$, however we can provide that $\epsilon\|\nabla\varrho\|_{2}$ will tend to zero, i.e. 
\begin{equation}\label{rhoepsilon}
\epsilon\|\nabla\varrho\|_{2}\leq \sqrt{\epsilon}C(\alpha,m_{2}),
\end{equation}
for some constant $C$ independent of $\epsilon$.\\ 
Now we would like obtain integrability of the pressure with the power $2$, as previously independently of $\epsilon$ and, if possible, of  $m_{2}$.\\
Therefore the choise of an appropriate test function seems to be obvious:
\begin{eqnarray*}
\Phi=\mathcal{B}\Big(P(\varrho)-\{P(\varrho)\}\Big),\quad in\ \Omega\\
\Phi=0\quad at\ \partial\Omega
\end{eqnarray*}
where $\mathcal{B}$ is the Bogovskii operator. By Lemma 3.17 from \cite{NS} and the Poincare inequality we have:
\begin{eqnarray}
\|\Phi\|_{\bar{p}}\leq c(p,\Omega)\|P(\varrho)\|_{p},\quad \|\nabla\Phi\|_{p}\leq c(p,\Omega)\|P(\varrho)\|_{p}\label{bogus}\\
0<p<\infty,\quad \bar{p}=\left\{\begin{array}{lcl}
\frac{2p}{2-p}&if&p<2\\
arbitrary\geq 1&if&p=2\\
\infty&if&p>2.
\end{array}\right.\nonumber
\end{eqnarray} 
From this testing,  the following identity appears:
\begin{multline*}
\calka P(\varrho)^{2}=\frac{1}{|\Omega|}\left(\calka P(\varrho)\right)^{2}+\alpha\calka(\varrho v - hg)\Phi+\mu\calka\nabla v:\nabla\Phi+(\mu+\nu)\calka\dyv v\ \dyv\Phi\\
-\calka K(\varrho)\varrho v\varotimes v:\nabla\Phi+\epsilon\calka\nabla\varrho\nabla v\Phi=\sum_{i=1}^{6}I_{i}.
\end{multline*}
Now each term will be estimated separately.\\
(i) By the estimate (\ref{osz1}) and the definition of $P$ the first one comes strightforward
\begin{equation*}
I_{1}=\frac{1}{|\Omega|}\left(\calka P(\varrho)\right)^{2}\leq\frac{1}{|\Omega|}\left(\calka\varrho^{\gamma}\right)^{2}\leq C.
\end{equation*}
(ii) The relation (\ref{bogus}) together with the estimate (\ref{osz1}) imply
\begin{multline*}
I_{2}=\alpha\calka (\varrho v-hg) \Phi\ dx\leq C\alpha\left(\|\varrho\|_{\gamma}^{1/2}\|\varrho v^{2}\|_{1}^{1/2}+\|h\|_{\gamma}^{1/2}\|hg^{2}\|_{1}^{1/2}\right)\|P(\varrho)\|_{2}\\
\leq C\alpha\|P(\varrho)\|_{2}.
\end{multline*}
(iii) We also have $\|\nabla \Phi\|_{2}\leq\|\varrho^{\gamma}\|_{2}$, thus
\begin{multline*}
I_{3}+I_{4}=\mu\calka\nabla v\nabla\Phi+(\mu+\nu)\calka\dyv v\dyv\Phi\leq C\|v\|_{2}\|P(\varrho)\|_{2}\\ \leq C\alpha^{1/2}\|P(\varrho)\|_{2}.
\end{multline*}
(iv) Since the modulus of $K$ is less than 1, the H$\mathrm{\ddot{o}}$lder's inequality and imbedding mentioned above lead to
\begin{equation*}
I_{5}=\calka K(\varrho)\varrho v\varotimes v:\nabla\Phi\leq C\|\varrho\|_{\gamma}\|v\|_{1,2}^{2}\|P(\varrho)\|_{2}\leq C\alpha\|P(\varrho)\|_{2}.
\end{equation*}
(v) Finally, epmloying the H$\mathrm{\ddot{o}}$lder's inequality we may get that
\begin{equation*}
I_{6}=\epsilon\calka\nabla\varrho\nabla v\Phi\leq\epsilon\|\nabla\varrho\|_{q}\|v\|_{1,2}\|P(\varrho)\|_{2},
\end{equation*}
for some $q>2$. To get the estimate for $\|\nabla\varrho\|_{q}$ we need to interpret the approximate continuity equation as a Neumann-boundary problem
	\begin{equation}\label{NBC}
	\begin{array}{c}
	-\epsilon\Delta\varrho=\dyv b\quad in\ \Omega\\
	\frac{\partial\varrho}{\partial n}=b\cdot n\quad at\  \partial\Omega,
	\end{array}
	\end{equation}
with the right hand side
\begin{equation*}
	b=\alpha \mathcal{B}(K(\varrho)h-\varrho)-K(\varrho)\varrho v. 
\end{equation*}
From the classical theory we know that if $\partial\Omega$ is smooth enough and if $b\in (L_{p}(\Omega))^{2}$, then there exists the unique $\varrho\in W^{1}_{p}(\Omega)$ satisfying  (\ref{NBC}) in the weak sence, such that $\calka \varrho dx=const$. Moreover
	\begin{equation}\label{reg0}
	\|\nabla\varrho\|_{p}\leq\frac{c(p,\Omega)}{\epsilon}\|b\|_{p}.
	\end{equation}
Now, in our case assume that $\gamma>q>2$ then the q-norm of $b$ may be estimated as
\begin{equation}\label{bogdan}
\|b\|_{q}\leq\alpha(\|\varrho\|_{q}+\|h\|_{q})+C\|\varrho\|_{\gamma}\|v\|_{1,2}\leq C_{1}\alpha+C_{2}\sqrt{\alpha},
\end{equation}
thus the observation (\ref{reg0}) yields the following
\begin{equation*}
I_{6}=\epsilon\calka\nabla\varrho\nabla v\Phi\leq (C_{1}\alpha^{3/2}+C_{2}\alpha)\|P(\varrho)\|_{2}.
\end{equation*}
Gathering the estimates terms $I_{i}$ for $i=1,\ldots,\ 6$ one can easily see that
\begin{equation}\label{2gamma}
\|P(\varrho)\|_{2}\leq C\alpha^{3/2},
\end{equation}
where the constant $C$ does not depend on $\epsilon$ nor $m_{2}$.\\
Now our aim will be to estimate the norm of $\nabla v$ in  $L_{q}(\Omega)$ for $q\geq 2$. For this puropse we will apply to the system (\ref{Lame}) 
the following Lemma (for the proof, see \cite{MP} Lemma 3.3.).
\begin{lemat}\label{wreg}
	Let $1<p<\infty,\ \Omega\in C^{2},\ F\in(M_{2p/(p+2)})^{*},\ \mu>0, 2\mu+3\nu>0$. Then there exists the unique $w\in M_{p}$, solution to (\ref{Lame}). Moreover
		\begin{equation*}
		\|w\|_{1,p}\leq C(p,\Omega)\|F\|_{(M_{p/(p-1)})^{*}}.
		\end{equation*}
	If $\Omega\in C^{l+2},\  F\in W^{l}_{p}(\Omega),\  l=0,1,\ldots$ then $w\in W^{l+1}_{p}(\Omega)$ and
		\begin{equation*}
		\|w\|_{l+2,p}\leq C(p,\Omega)\|F\|_{l,p}.
		\end{equation*}
	\end{lemat}	
If we consider the approximate momentum equation as a part of Lame system with $w=v$ we will get the estimate for the norm of $\nabla v$ in $L_{q}(\Omega)$
\begin{multline*}
\|\nabla v\|_{q}\leq C(\alpha\|\varrho v\|_{2q/(q+2)}+\alpha\|hg\|_{2q/(q+2)})+\|K(\varrho)\varrho v\varotimes v\|_{q}+\|P(\varrho)\|_{q}\\+\epsilon\|\nabla\varrho\nabla v\|_{2q/(q+2)}).
\end{multline*}
Recalling $\gamma>2$, by (\ref{osz1}) and by (\ref{vH1}) one gets\\
\begin{equation*}
\alpha\|\varrho v\|_{2q/(q+2)}+\alpha\|hg\|_{2q/(q+2)}\leq C \alpha(\|\varrho v\|_{2}+\|hg\|_{2})\leq C\alpha^{3/2}.
\end{equation*}
By the definition of $P$ and the H$\mathrm{\ddot{o}}$lder's inequality we also have
\begin{equation*}
\|K(\varrho)\varrho v\varotimes v\|_{q}\leq C\|P(\varrho)\|_{q/\gamma}^{\gamma}\|v\|^{2}_{1,2}\leq C\alpha\|P(\varrho)\|^{1/\gamma}_{q/\gamma}.
\end{equation*}
At this step there is a need to include the estimates depending on the parameter $m_{2}$, more precisely we will use 
\begin{eqnarray*}
&{}&\|P(\varrho)\|_{q}\leq \|P(\varrho)\|_{\infty}^{1-2/q}\|P(\varrho)\|_{2}^{2/q}\leq C \alpha^{3/q}m_{2}^{(1-2/q)\gamma},\\
&{}&\epsilon\|\nabla\varrho\nabla v\|_{2q/(q+2)}\leq\epsilon\|\nabla\varrho\|_{q}\|v\|_{1,2}\leq C\alpha^{3/2}m_{2},
\end{eqnarray*}
where the last inequality is obtained by replacing in (\ref{bogdan}) the norms of $\varrho$ by $\|\varrho\|_{\infty}\leq m_{2}$ if $q\geq\gamma$; for $q<\gamma$ we can use the estimate for the $L_{\gamma}$-norm of $\varrho$.\\
Summarising, we have shown that $\|\nabla v\|_{q}\leq C(m_{2},\alpha)$ with a constant $C(m_{2},\alpha)$ independent of $\epsilon$. Particulary for $2<q<\gamma$ one has
\begin{equation}\label{gradq}
\|\nabla v\|_{q}\leq C(\alpha+\alpha^{3/2}+\alpha^{3/q}m_{2}^{(1-2/q)\gamma}).
\end{equation}
Before  passing to the zero limit with $\epsilon$ we will compute {\it a priori} estimate of the vorticity
\begin{equation*}
\omega=\rot v=\frac{\partial v_{2}}{\partial x_{1}}-\frac{\partial v_{1}}{\partial x_{2}}.
\end{equation*}
Differentiating $n\cdot v=0$ at $\partial \Omega$ with respect to the length parameter and combining it with the last boundary condition in the system (\ref{approx}) we obtain:
\begin{equation*}
\omega=\left(2\chi-\frac{f}{\mu}\right)v\cdot\tau\quad at\ \partial\Omega.
\end{equation*} 
Taking the rotation of $(\ref{approx})_{2}$, we get
\begin{equation}
-\mu\Delta\omega=-\alpha\rot(hg-\varrho v)-\rot\dyv(K(\varrho)\varrho v\varotimes v)-\epsilon\rot(\nabla\varrho\nabla v).\label{rot}
\end{equation}
Denote 
\begin{equation}\label{omegi}
\omega=\omega_{1}+\omega_{2},
\end{equation}
where $\omega_{1},\ \omega_{2}$ satisfy:
\begin{eqnarray*}
&{}& -\mu\Delta\omega_{1}=-\rot\dyv(K(\varrho)\varrho v\varotimes v)\quad in\ \Omega,\\
&{}& \omega_{1}=0\quad at\ \partial\Omega,\\
&{}& -\mu\Delta\omega_{2}=-\alpha\rot(hg-\varrho v)-\epsilon\rot(\nabla\varrho\nabla v)\quad in\ \Omega,\\
&{}& \omega_{2}=\left(2\chi-\frac{f}{\mu}\right)v\cdot\tau\quad at\ \partial\Omega.
\end{eqnarray*}
For the weak solutions $\omega_{1},\omega_{2}$ of the above problems one get the following estimates:
\begin{equation*}
\|\omega_{1}\|_{q}\leq C\|K(\varrho)\varrho v\varotimes v\|_{q}\leq C\alpha
\end{equation*}
where for $q<\gamma$, $C$ is independent of $m_{2}$ and for $q>\gamma$,  $C=C_{0}m_{2}^{1-\gamma/q}$,
\begin{equation*}
\|\omega_{2}\|_{1,q}\leq C(\alpha \|hg\|_{q}+\alpha\|\varrho v\|_{q}+\epsilon\|\nabla\varrho\nabla v\|_{q})+C(\Omega)\|v\cdot\tau\|_{1-1/q,q,\partial\Omega},
\end{equation*}
thus for $q< 2$, the H$\mathrm{\ddot{o}}$lder's inequality, the imbedding $W^{1/2}_{2}(\partial\Omega)\subset W^{1-1/q}_{q}(\partial\Omega)$ and the trace theorem imply
\begin{multline*}
\|\omega_{2}\|_{1,q}\leq C(\alpha \|hg\|_{\frac{2\gamma}{\gamma+1}}+\alpha\|\varrho v\|_{\frac{2\gamma}{\gamma+1}}+\epsilon\|\nabla\varrho\|_{2q/(2-q)}\|\nabla v\|_{2})+C(\Omega)\|v\|_{1,2}\\\leq C(\alpha+\alpha^{2}m_{2})+C(\Omega)\alpha^{1/2},
\end{multline*}
for $q\geq 2$ we must use $m_{2}$-dependent estimates of gradient of $v$ in higher norms, thus
\begin{equation*}
\|\omega_{2}\|_{1,q}\leq C(\alpha,m_{2})
\end{equation*}
and the dependence of $m_{2}$ is higher then linear.
\section{Passage to the limit}
This section is devoted to the passage with $\epsilon\rightarrow 0$ in the system (\ref{approx}). Recall that so far we have obtained the following estimates:
\begin{eqnarray}
&{}&\|\varrho_{\epsilon}\|_{\infty}\leq m_{2}, \quad\|v_{\epsilon}\|_{1,2}\leq C\alpha,\\
&{}&\|P(\varrho_{\epsilon})\|_{2}\leq C\alpha^{3/2}\label{normaP}\\
&{}&\|v_{\epsilon}\|_{1,q}+\epsilon^{1/2}\|\nabla\varrho_{\epsilon}\|_{2}\leq C(m_{2},\alpha)\quad for\ 1\leq q<\infty,\\
&{}&\epsilon\|\nabla\varrho_{\epsilon}\nabla v_{\epsilon}\|_{q}\leq C(m_{2},\alpha)\quad for\ 1\leq q<\infty.
\end{eqnarray} 
The two last estimates together with the interpolation inequality imply that for $\delta$ sufficiently small we additionally have:
\begin{equation*}
\epsilon^{1-\delta}\|\nabla\varrho_{\epsilon}\nabla v_{\epsilon}\|_{q}\leq C(m_{2},\alpha)\quad for\ 1\leq q<\infty.
\end{equation*}
Therefore, at least for an appropriately chosen subsequence:
\begin{eqnarray*}
&{}&\varrho_{\epsilon}\rightharpoonup^{*}\varrho\quad in\ L_{\infty}(\Omega),\\
&{}&P(\varrho_{\epsilon})\rightharpoonup\overline{P(\varrho)}\quad in\ L_{2}(\Omega),\\
&{}&v_{\epsilon}\rightharpoonup v\quad in\ W^{1}_{q}(\Omega),\\
&{}&\epsilon\nabla\varrho_{\epsilon}\rightarrow 0\quad in\ L_{2}(\Omega),\\
&{}&\epsilon\nabla\varrho_{\epsilon}\nabla v_{\epsilon}\rightarrow 0\quad in\ L_{q}(\Omega)\ for\ 1\leq q<\infty,
\end{eqnarray*}
where the line over a term denotes its weak limit.\\
These information allow us to pass to the limit in our approximative system:
\begin{equation}\label{granica}
\begin{array}{r}
\alpha\left(\varrho-\overline{hK(\varrho)}\right)+\dyv(\overline{K(\varrho)\varrho} v)=0\\
\alpha\left(\varrho v-hg\right) +\dyv(\overline{K(\varrho)\varrho} v\varotimes v)-\mu\Delta v-(\mu+\nu)\nabla\dyv v+\nabla \overline{P(\varrho)}=0\\
v\cdot n=0\quad at\quad \partial\Omega,\\
n\cdotp T(v,\overline{P(\varrho)})\cdotp\tau+fv\cdotp\tau=0\quad at\quad\ \partial\Omega.
\end{array}
\end{equation}
To show that we have realy found the solution to our initial problem there left several questions that need to find the answer. \\
Firstly, if we can get rid of $K(\varrho)$ that remains at several places, i.e. if we can prove that $K(\varrho)= 1$ a.e. in $\Omega$. This, as we shall see below, is equivalent with showing that there can be suitably chosen constant $m$ sufficiently large but still sharply smaller than the {it a priori} bound for a density, such that  the measure of the set
\begin{equation*}
\{x\in\Omega:\varrho_{\epsilon_{n}}(x)>m\}
\end{equation*}
tends to zero for some subsequence $\epsilon_{n}\rightarrow 0^{+}$. Indeed, as this implies that for any smooth function $\eta$ one get
\begin{equation*}
\calka\varrho_{\epsilon_{n}}K(\varrho_{\epsilon_{n}})\eta\ dx=\calka\varrho_{\epsilon_{n}}\eta\ dx+\int_{\{\varrho_{\epsilon_{n}}>m_{1}\}}(K(\varrho_{\epsilon_{n}})-1)\varrho_{\epsilon_{n}}\eta\ dx.
\end{equation*}
If we choose $m_{1}$ sufficiently close to $m_{2}$ and additionally assure that $m<m_{1}$ then the last term on the right hand side disappears as $\epsilon_{n}$ goes to 0, and thus we truly have
\begin{equation*}
\lim_{\epsilon_{n}\rightarrow 0^{+}}\calka\varrho_{\epsilon_{n}}K(\varrho_{\epsilon_{n}})\eta\ dx=\calka\varrho\eta\ dx,\quad \forall\eta\in C^{\infty}(\Omega).
\end{equation*}
The next difficulty concerns the convergence in the nonlinear term i.e. is it true that $\overline{P(\varrho)}=P(\varrho)$. The positive answer can be obtained in a rather standard way, and at the stage when one already knows that $K(\varrho)=1$ it reduces to proving the strong convergence for the density sequence.\\
Finally, what does the condition $(\ref{granica})_{4}$ mean, in other words in which sense is it satisfied?
Having solved the two previous problem this is quite easy to see that this boundary condition can be recovered while passing to the limit in a weak formulation corresponding to the momentum equation.\\
Now our aim will be to precisely justify the considerations developed above. For this purpose we will adapt a kind of technique widely used for these type of problems, more precisely we will take advantage of some properties of the effective viscous flux denoted in this paper by $G$.\\
Introducing the Helmholtz decomposition  of the velocity vector field defined as:
\begin{equation}\label{decomposition}
v=\nabla\phi+\nabla^{\bot} A,
\end{equation}
where the divergence-free part $\nabla^{\bot} A=\left(-\frac{\partial}{\partial x_{2}},\frac{\partial}{\partial x_{1}}\right)A$ and the gradient part $\phi$ are given by:
\begin{equation}\label{rozklad}
\left\{\begin{array}{ll}
\Delta A=\rot v\quad &in\ \Omega\\
\nabla^{\bot} A\cdot =0\quad &at\ \partial\Omega
\end{array}\right. ,\quad
\left\{\begin{array}{ll}
\Delta \phi=\dyv v\quad &in\ \Omega,\\
\frac{\partial\phi}{\partial n}=0\quad &at\ \partial\Omega
\end{array}\right. ,
\end{equation}
we can transform the limit equation $(\ref{granica})_{2}$ into the form:
\begin{equation}\label{definG}
\nabla G=\alpha hg-\alpha\varrho v-\dyv(\overline{K(\varrho)\varrho}v\varotimes v)+\mu\Delta\nabla^{\bot}A,
\end{equation}
where $\nabla G=\nabla\left(-(2\mu+\nu)\Delta\phi+\overline{P(\varrho)}\right)$. 
By the observation $\calka G dx=\calka\overline{P(\varrho)}dx\leq\infty$, we control the mean value of $G$ and thus the expression
\begin{equation*}
G=(2\mu+\nu)\Delta\phi+\overline{P(\varrho)}
\end{equation*}
may be accepted as a correct definition of $G$.\\
Due to (\ref{rozklad}), and the classical theory for the laplacian supplemented by the Neuman-boundary condition 
\begin{equation*}
\|G\|_{2}\leq C(\|\nabla v\|_{2}+\|\overline{P(\varrho)}\|_{2})\leq C(\alpha).
\end{equation*}
The next goal is to show the boudedness of the $L_{\infty}$ norm of $G$. By the fact that the mean value of $G$ is controlled we can  employ the Poincare's inequality and the Sobolev embedding theorem it is sufficient to prove the following fact:
\begin{lemat}
For $q>2$ we have:
\begin{equation}
\|\nabla G\|_{q}\leq C(\alpha,m_{2}).
\end{equation}
\end{lemat}
\dow
By virtue of (\ref{definG}) 
\begin{equation}\label{Gdelta}
\|\nabla G\|_{q}\leq C \alpha\|hg\|_{q}+\alpha\|\varrho v\|_{q}+\|\dyv(\overline{K(\varrho)\varrho}v\varotimes v)\|_{q}+\mu\|\Delta\nabla^{\bot}A\|_{q}.
\end{equation}
For $q<\gamma$ may we certainly write that
\begin{equation*}
\alpha\|hg\|_{q}+\alpha\|\varrho v\|_{q}\leq C\alpha\|v\|_{1,2}\leq C\alpha^{3/2},
\end{equation*}
by the continuity equation, the estimates (\ref{osz2}) and (\ref{nonnegativ}) we get
\begin{multline*}
\|\dyv(\overline{K(\varrho)\varrho}v\varotimes v)\|_{q}\leq\|\overline{K(\varrho)\varrho}v\cdot\nabla v\|_{q}+\alpha\|\overline{hK(\varrho)}v\|_{q}+\alpha\|\varrho v\|_{q}\\
\leq Cm_{2}\|\nabla v\|_{q}^{2} +C\alpha^{3/2},
\end{multline*}
thus, the estimate (\ref{gradq}) of $\|\nabla v\|_{q}$ for $\gamma>q>2$ leads to
\begin{equation*}
\|\dyv(\overline{K(\varrho)\varrho}v\varotimes v)\|_{q}\leq C(\alpha^{3/2}+\alpha^{2}+\alpha^{3}+\alpha^{6/q}m_{2}^{1+2(1-2/q).\gamma})
\end{equation*}
The last term in (\ref{Gdelta}) is bounded by the same constant, since
\begin{multline*}
\|\Delta\nabla^{\bot}A\|_{q}\leq\|\nabla\omega\|_{q}\leq\alpha\|hg\|_{q}+\alpha\|\varrho v\|_{q}+\|\dyv(\overline{K(\varrho)\varrho}v\varotimes v)\|_{q}+\\+C\|v\cdot \tau\|_{1-{1/q},2+\delta,\partial\Omega},
\end{multline*}
where $\omega$ is a weak solution to (\ref{rot}) with a corresponding boundary condition after passing with $\epsilon$ to $0$, i.e. it satisfies
\begin{eqnarray*}
-\mu\Delta\omega&=&-\alpha\rot(hg-\varrho v)-\rot\dyv(\overline{K(\varrho)\varrho} v\varotimes v)\quad in\ \Omega\\
\omega&=&\left(2\chi-\frac{f}{\mu}\right)v\cdot\tau\quad at\ \partial\Omega.
\end{eqnarray*}
\begin{flushright}
$\Box$
\end{flushright}
For $q$  such small that $\gamma>1+2(1-2/q)\gamma$ we have then proved that
\begin{equation}\label{Ginf}
\|G\|_{\infty}\leq C(\alpha)m_{2}^{\gamma-\delta},
\end{equation}
with $\delta=\gamma(4/q-1)-1>0$ and $C(\alpha)=6/q$\\
We will now apply the analogical decomposition for the approximative system (\ref{approx}), i.e.
\begin{equation*}
v_{\epsilon}=\nabla\phi_{\epsilon}+\nabla^{\bot} A_{\epsilon}.
\end{equation*}
Similarly as previously this leads to relation
\begin{multline}
\nabla G_{\epsilon}=(2\mu+\nu)\Delta\phi_{\epsilon}+P(\varrho_{\epsilon})\\=\alpha hg-\alpha\varrho_{\epsilon} v_{\epsilon}-\dyv(K(\varrho_{\epsilon})\varrho_{\epsilon}v_{\epsilon}\varotimes v_{\epsilon})-\epsilon\nabla\varrho_{\epsilon}\nabla v_{\epsilon}+\mu\Delta\nabla^{\bot}A_{\epsilon}.
\end{multline}
We are then able to prove that if $\epsilon\rightarrow 0^{+}$ the following lemma holds
\begin{lemat}\label{strong}
$G_{\epsilon}\rightarrow G$ strongly in $L_{2}(\Omega)$.
\end{lemat}
\dow
We will use the following fact:
\begin{equation*}
If\quad \nabla (G_{\epsilon}-G)\rightharpoonup 0\quad in\ L_{2},\quad then\quad G_{\epsilon}-G\rightarrow const\quad in\ L_{2}, 
\end{equation*}
and next we can show that at least for some subsequence $\epsilon_{n}\rightarrow 0$ the constant is indeed equal zero
\begin{equation*}
\calka(G_{\epsilon}-G)=\calka\Delta(\phi-\phi_{\epsilon})\rightarrow 0
\end{equation*}
since $\frac{\partial \phi}{\partial n}=\frac{\partial \phi_{\epsilon}}{\partial n}=0$ at $\partial\Omega$.\\

This allows us to focus on showing the weak convergence, we have
\begin{multline}\label{Gzbieg}
\nabla(G_{\epsilon}-G)=\mu\Delta\nabla^{\bot}(A^{\epsilon}-A)-\alpha(\varrho_{\epsilon}v_{\epsilon}-\varrho v)\\-(\dyv(K(\varrho_{\epsilon})\varrho_{\epsilon}v_{\epsilon}\varotimes v_{\epsilon})-\dyv(\overline{K(\varrho)\varrho}v\varotimes v))-\epsilon\nabla\varrho_{\epsilon}\nabla v_{\epsilon}.
\end{multline}
The second term on the right hand side converges to 0 weakly in $L_{2}$ owing to the strong convergence of $v_{\epsilon}\rightarrow v$ in $L_{q}$ for any $0\leq q\leq \infty$ and by the boundedness of $\varrho_{\epsilon}$ in $L_{\infty}$.\\
The last term converges to zero even strongly in $L_{2}$. Now, by the continuity equation, the third term may be written in the form
\begin{multline*}
\dyv(K(\varrho_{\epsilon})\varrho_{\epsilon}v_{\epsilon}\varotimes v_{\epsilon})-\dyv(\overline{K(\varrho)\varrho}v\varotimes v)=\alpha hK(\varrho_{\epsilon})v_{\epsilon}-\varrho_{\epsilon} v_{\epsilon}+\epsilon\Delta\varrho_{\epsilon}v_{\epsilon}\\+\alpha\varrho v-\alpha\overline{hK(\varrho)}v+K(\varrho_{\epsilon})\varrho_{\epsilon}v_{\epsilon}\cdot\nabla  v_{\epsilon}-\overline{K(\varrho)\varrho}v\cdot \nabla v,
\end{multline*}
due to the argument explained above we need to justify the convergence only for two terms. Firstly note that $\epsilon\Delta\varrho_{\epsilon}v_{\epsilon}$ converges to $0$ strongly in $W^{-1}_{2}(\Omega)$. Secondly, since $\nabla(v_{\epsilon}-v)\rightharpoonup 0$ weakly in $L_{2}(\Omega)$ we obtain the same information for $K(\varrho_{\epsilon})\varrho_{\epsilon}v_{\epsilon}\cdot\nabla  v_{\epsilon}-\overline{K(\varrho)\varrho}v\cdot \nabla v$.\\
In order to substantiate, taht the first term in (\ref{Gzbieg}) also tends to $0$ we observe that 
\begin{equation}\label{Aomega}
\Delta\nabla^{\bot}(A^{\epsilon}-A)=\nabla^{\bot}(\omega_{\epsilon}-\omega),
\end{equation} 
and that the function $\omega_{\epsilon}-\omega$ satisfies the system of equations
\begin{eqnarray*}
-\mu\Delta(\omega_{\epsilon}-\omega)&=&-\alpha\rot(\varrho_{\epsilon} v-\varrho v)-\rot\dyv(K(\varrho_{\epsilon})\varrho_{\epsilon} v_{\epsilon}\varotimes v_{\epsilon}- \overline{K(\varrho)\varrho} v\varotimes v)\\
&{}&-\epsilon\ \rot(\nabla\varrho_{\epsilon}\nabla v_{\epsilon})\quad in\ \Omega\\
\omega_{\epsilon}-\omega&=&\left(2\chi-\frac{f}{\mu}\right)(v_{\epsilon}-v)\cdot\tau\quad at\ \partial\Omega.
\end{eqnarray*}
Repeating the same reasoning as in case of (\ref{omegi}) and above explications  we can show, that $\nabla(\omega_{\epsilon}-\omega)$ consists of two parts. One of them converges to $0$ strongly in $W^{-1}_{2}(\Omega)$ and the other converges weakly in $L_{2}(\Omega)$. Thus, by (\ref{Aomega}), we get the same for $\Delta\nabla^{\bot}(A^{\epsilon}-A)$ and therefore the proof of lemma is complete.
\begin{flushright}
$\Box$
\end{flushright}
Provided with these information we can show the final argument for $K(\varrho)$ to be equal $1$
\begin{lemat}\label{lematK}
Let $\kappa>0$ and let $m$ satisfy
\begin{equation*}
\|G\|_{\infty}^{1/\gamma}<m<m_{1}\quad and\quad \frac{m^{\gamma+1}}{m_{2}}-\|G\|_{\infty}-\alpha(2\mu+\nu)\geq\kappa>0
\end{equation*}
then we have
\begin{equation*}
\lim_{\epsilon_{n}\rightarrow 0^{+}}|\{x\in\Omega:\varrho_{\epsilon_{n}}(x)>m\}|=0.
\end{equation*}
\end{lemat}
\dow\\
The difference with respect to the Lemma 4.3 from \cite{MP} is that the rate of convergence here clearly must depend on $\alpha$ and thus we pass with $\epsilon$ to 0 when $\alpha$ is set.\\
First observe that the assumptions of our lemma are satisfied. Indeed, as the difference $m_{2}-\|G\|^{1/\gamma}_{\infty}$ increases with $m_{2}$.
Next, we introduce a function $M(\cdot)\in C^{1}(\mathbb{R})$ given by
\begin{equation*}
M(\varrho)=\left\{
\begin{array}{ll}
1 &\quad \varrho\leq m,\\
0 &\quad \varrho\geq m+1,\\
\in(0,1) &\quad \varrho\in(m,m+1),
\end{array}
\right.
\end{equation*}
where $M'(\varrho)<0$ in $(m,m+1)$ and $m+1<m_{1}$.\\
We multiply the approximate continuity equation by $M^{l}(\varrho_{\epsilon})$ for some $l\in\mathbb{N}$ and we observe
\begin{multline}
\alpha\calka M^{l}(\varrho_{\epsilon})\left(\varrho-hK(\varrho)\right)dx+\calka M^{l}(\varrho_{\epsilon})\dyv(K(\varrho)\varrho v)dx=\epsilon\calka M^{l}(\varrho_{\epsilon})\Delta\varrho\ dx\\
=-\epsilon l\calka M'(\varrho_{\epsilon})M^{l-1}(\varrho_{\epsilon})|\nabla\varrho_{\epsilon}|^{2}\ dx\geq 0.
\end{multline}
By integrating the second term on the left hand side by parts twice  (the boundary terms disappear due to the definition of $M(\cdot)$) one gets
\begin{multline*}
\calka\left(\int_{0}^{\varrho_{\epsilon}}tM^{l-1}(t)M'(t)dt\right)\dyv v_{\epsilon}\ dx\\ \geq\frac{\alpha}{l}\calka\left(hK(\varrho_{\epsilon})-\varrho_{\epsilon}\right)\ dx+\frac{\alpha}{l}\calka\left(\varrho_{\epsilon}-hK(\varrho_{\epsilon})\right)\left(1-M^{l}(\varrho_{\epsilon})\right)\ dx.
\end{multline*}
The first therm on the right hand side cancels due to the Theorem \ref{twr:approx}. We can replace $\dyv v_{\epsilon}$ according to the definition of $G_{\epsilon}$, then we have
\begin{multline*}
\calka\left(\int_{0}^{\varrho_{\epsilon}}tM^{l-1}(t)M'(t)dt\right) \left(G_{\epsilon}-P(\varrho_{\epsilon})\right)\ dx\\ \leq-\frac{\alpha(2\mu+\nu)}{l}\calka\left(\varrho_{\epsilon}-hK(\varrho_{\epsilon})\right)\left(1-M^{l}(\varrho_{\epsilon})\right)\ dx.
\end{multline*}
Since $M'(t)$ is negative, supported in $(m,m+1)$ and $m+1<m_{1}\rightarrow m_{2}^{-}$ the following inequality holds true
\begin{multline*}
-m\calka\left(\int_{0}^{\varrho_{\epsilon}}M^{l-1}(t)M'(t)dt\right)P(\varrho_{\epsilon})\ dx\\ \leq m_{2}\calka\left|-\int_{0}^{\varrho_{\epsilon}}M^{l-1}(t)M'(t)dt\right||G_{\epsilon}|\ dx+\frac{\alpha(2\mu+\nu)}{l}\calka\left|\varrho_{\epsilon}-hK(\varrho_{\epsilon})\right|\left(1-M^{l}(\varrho_{\epsilon})\right)\ dx.
\end{multline*}
After integration of the internal term we claim to conclusion that the above expression is different then $0$ only for a subset of $\Omega$, $\{\varrho_{\epsilon}>m\}$, thus
\begin{multline}\label{cos}
\frac{m}{m_{2}}\int_{\{\varrho_{\epsilon}>m\}}(1-M^{l}(\varrho_{\epsilon}))P(\varrho_{\epsilon})\ dx\\
\leq\int_{\{\varrho_{\epsilon}>m\}}(1-M^{l}(\varrho_{\epsilon}))|G_{\epsilon}|\ dx +\frac{\alpha(2\mu+\nu)}{m_{2}}\int_{\{\varrho_{\epsilon}>m\}}\left|\varrho_{\epsilon}-hK(\varrho_{\epsilon})\right|\left(1-M^{l}(\varrho_{\epsilon})\right)\ dx.
\end{multline}
Now for each $\delta>0$ we can find such sufficiently large number $l\in \mathbb{N}$, $l=l(\delta, \epsilon)$ that
\begin{equation}\label{costam}
\|M^{l}(\varrho_{\epsilon})\|_{L_{2}(\{\varrho_{\epsilon}>m\})}\leq \delta,
\end{equation}
since $M(\varrho_{\epsilon})$ is less then $1$ for $\varrho_{\epsilon}>m$. This allows us to rewrite the inequality (\ref{cos}) in the following form
\begin{multline*}
\frac{m^{\gamma+1}}{m_{2}}\left|\{\varrho_{\epsilon}>m\}\right|\leq\frac{m}{m_{2}}\|M^{l}(\varrho_{\epsilon})\|_{L_{2}(\{\varrho_{\epsilon}>m\})}\|P(\varrho_{\epsilon})\|_{L_{2}(\{\varrho_{\epsilon}>m\})}\\+ C(|\Omega|)\|G-G_{\epsilon}\|_{2}+\|G\|_{\infty} \left|(\{\varrho_{\epsilon}>m\}\right|+\alpha(2\mu+\nu)\left|(\{\varrho_{\epsilon}>m\}\right|,
\end{multline*}
where the term on the left is a consequence of the definition of $P(\cdot)$ and the limits of integration. By (\ref{costam}) and the bound (\ref{normaP}) we therefore may write
\begin{multline*}
\left(\frac{m^{\gamma+1}}{m_{2}}-\|G\|_{\infty} -\alpha(2\mu+\nu)\right)\left|\{\varrho_{\epsilon}>m\}\right|\leq\frac{Cm}{m_{2}}\delta\alpha^{3/2}+ C(|\Omega|)\|G-G_{\epsilon}\|_{2}.
\end{multline*}
Under our assumptions, the expression in the brackets is separated from $0$ and at least for a suitably chosen subsequence $\epsilon_{n}\rightarrow 0^{+}$ Lemma \ref{strong} guarantees that
\begin{equation*}
\lim_{\epsilon_{n}\rightarrow 0^{+}}\left|\{\varrho_{\epsilon_{n}}>m\}\right|\leq\frac{Cm}{m_{2}}\delta\alpha^{3/2}.
\end{equation*}
As $\delta$ may be arbitrary small and $\alpha=const$, we truly have
\begin{equation*}
\lim_{\epsilon_{n}\rightarrow 0^{+}}\left|\{\varrho_{\epsilon_{n}}>m\}\right|=0.
\end{equation*}
\begin{flushright}
$\Box$
\end{flushright}
This fact, as it was already mentioned before, completes  justification that $K(\varrho)=1$ a.e. in $\Omega$.\\
The second problem to solve was to show that  $\overline{P(\varrho)}=P(\varrho)$. For this purpose we multiply the approximate continuity equation by the function $\ln\frac{m_{2}}{\varrho_{\epsilon}+\delta}$ for $\delta>0$ and integrate over $\Omega$. Like in the proof of last lemma, we observe
\begin{multline}
\alpha\calka \ln\frac{m_{2}}{\varrho_{\epsilon}+\delta}\left(\varrho-h\right)dx+\calka \ln\frac{m_{2}}{\varrho_{\epsilon}+\delta}\dyv(\varrho v)dx\\=\epsilon\calka \ln\frac{m_{2}}{\varrho_{\epsilon}+\delta}\Delta\varrho\ dx
=\epsilon l\calka \frac{|\nabla\varrho_{\epsilon}|^{2}}{\varrho_{\epsilon}+\delta}\ dx\geq 0.
\end{multline}
Similarly as previously we integrate by parts, pass with $\delta\rightarrow 0^{+}$, substitute $G_{\epsilon}$ from the definition and pass with $\epsilon\rightarrow 0^{+}$ to get
\begin{equation}\label{Ggora}
\calka\overline{P(\varrho)\varrho}\ dx+(2\mu+\nu)\alpha\calka\overline{(\varrho-h)\ln\varrho}\ dx\leq\calka G\varrho\  dx.
\end{equation}
From now on we will seek to reverse the sign of above inequality. We will use the fact that the limit continuity equation works with any smooth function up to the boundary. To indicate an appropriate one we first introduce  the distribution:
\begin{equation*}
v\cdot\nabla\varrho=\dyv(\varrho v)-\varrho\dyv v.
\end{equation*} 
Then let us recall the following lemma (for the proof consult \cite{NNP}).
\begin{lemat}\label{8,11}
Let $\Omega\in C^{0,1},\ v\in W^{1}_{q}(\Omega),\  \varrho\in L_{p}(\Omega), 1<p,q<\infty$, $v\cdot\nabla\varrho\in L_{s}(\Omega)$, $1/s=1/p+1/q$. Then there exists $\varrho_{n}\in C^{\infty}(\overline{\Omega})$ such that
\begin{equation*}
v\cdot\nabla\varrho_{n}\rightarrow v\cdot\nabla\varrho\ in\ L_{s}(\Omega)\quad and\quad\varrho_{n}\rightarrow\varrho\ in\ L_{p}(\Omega).
\end{equation*}
\end{lemat}
\noindent For such a $\varrho_{n}$ one gets
\begin{equation*}
\calka\dyv(\varrho_{n}v)dx=\calkad\varrho_{n}v\cdot ndS=0,
\end{equation*}
thus passing with $n\rightarrow\infty$ our lemma provides that
\begin{equation*}
\calka\varrho\dyv v dx=-\calka v\cdot\nabla\varrho dx.
\end{equation*}
Note that a function $\ln\frac{\delta}{\varrho_{n}+\delta}$ for $\delta>0$ is an admissible test function as it follows from the proof of Lemma \ref{8,11} that $0\leq\varrho_{n}\leq m_{2}$, hence we get
\begin{equation*}
\alpha\calka(h-\varrho)\ln\frac{\delta}{\varrho_{n}+\delta}=\calka\varrho v\frac{\nabla\varrho_{n}}{\varrho_{n}+\delta}.
\end{equation*}
 We may now pass with $n\rightarrow\infty$
 \begin{equation*}
 \alpha\calka(h-\varrho)\ln\frac{\delta}{\varrho+\delta}=\calka\frac{\varrho v\cdot\nabla\varrho}{\varrho+\delta}.
 \end{equation*}
Next we also want to pass with $\delta\rightarrow 0^{+}$, since $\calka(\varrho-h)\ln\delta\ dx=0$, the only difficult term is $\alpha\calka h\ln(\varrho+\delta)$, but it can be solved by the Lebesgue monotone convergence theorem, then we obtain
 \begin{equation*}
 \alpha\calka h\ln\varrho=\alpha\calka \varrho\ln\varrho-\calka v\cdot\nabla\varrho=\alpha\calka \varrho\ln\varrho+\calka \varrho\dyv v.
 \end{equation*}
 Finally, recalling the definition of $G$ one gets
 \begin{equation}\label{Gdol}
 \calka G\varrho\ dx=(2\mu+\nu)\alpha\calka (\varrho-h)\ln\varrho\  dx +\calka \overline{P(\varrho)}\varrho\ dx.
 \end{equation}
 The information contained in (\ref{Ggora}), (\ref{Gdol}) together imply
 \begin{equation}\label{Gtot}
 \calka\overline{P(\varrho)\varrho}\ dx+(2\mu+\nu)\alpha\calka\overline{(\varrho-h)\ln\varrho}\ dx\leq(2\mu+\nu)\alpha\calka (\varrho-h)\ln\varrho\  dx +\calka \overline{P(\varrho)}\varrho\ dx.
 \end{equation}
 The convexity of functions $\varrho\ln(\varrho)$  and $-h\ln(\varrho)$ ensure lower semicontinuity of the functional $\calka(\varrho-h)\ln(\varrho)\  dx$, in other words
 \begin{equation}
 \calka (\varrho-h)\ln\varrho\  dx\leq\calka\overline{(\varrho-h)\ln\varrho}\ dx.
 \end{equation}
Therefore  (\ref{Gtot}) reduces to
\begin{equation}\label{pompka}
\calka\overline{P(\varrho)\varrho}\ dx\leq\calka \overline{P(\varrho)}\varrho\ dx.
\end{equation}
By the 'standard arguments' we show that $\varrho\overline{\varrho^{\gamma}}\leq\overline{\varrho^{\gamma+1}}$ which together with (\ref{pompka}) yield
\begin{equation}
\overline{\varrho^{\gamma}}\varrho=\overline{\varrho^{\gamma+1}}.
\end{equation}
Next, we may also show that $\overline{\varrho^{\gamma}}^{(\gamma+1)/\gamma}(x)\leq\overline{\varrho^{\gamma+1}}(x)$ and $\varrho(x)\leq\overline{\varrho^{\gamma}}^{1/\gamma}$ for a.a. $x\in\Omega$ which easily imply that
\begin{equation}
\varrho^{\gamma}(x)=\overline{\varrho^{\gamma}}(x).
\end{equation}
Since $L_{\gamma}(\Omega)$ is a uniformly convex Banach Space for $\gamma>1$, $\varrho_{\epsilon}\rightharpoonup \varrho$ weakly in $L_{\gamma}(\Omega)$ and $\|\varrho_{\epsilon}\|_{\gamma}^{\gamma}\rightarrow\|\varrho\|_{\gamma}^{\gamma}$ we may deduce, that $\varrho_{\epsilon}\rightarrow\varrho$ strongly in $L_{\gamma}(\Omega)$. Thus in turn implies, that for some subsequence $\varrho_{\epsilon}\rightarrow\varrho$ a.e. in $\Omega$ and the condition $\|\varrho_{\epsilon}\|_{L_{\infty}(\Omega)}$ guarantees the uniform integrability of the sequence $\{\varrho_{\epsilon_{n}}\}_{n=1}^{\infty}$ which together with the Vitali's convergence theorem leads to the strong convergence of the approximate densities to the function $\varrho$ in  $L_{p}(\Omega)$ for any $1\leq p<\infty$.\\
\begin{remark}
The density obtained in the above procedure is bounded by $m$ as we could see in lemma \ref{lematK}. Now, by taking $\kappa$ sufficiently small and $m_{1}, m_{2}$ sufficiently close to $m$ the estimate (\ref{Ginf}) for $q\rightarrow 2^{+}$ with condition imposed on the assumptions of Lemma \ref{lematK} will imply that
\begin{equation*}
\|\varrho\|_{\infty}\leq C(\alpha)^{3/(\gamma-1)},
\end{equation*}
in particular, $\|\varrho\|_{\infty}\leq C(\triangle t)^{-3}$.
\end{remark}
Theorem \ref{theorem1} is proved.
\begin{flushright}
$\Box$
\end{flushright}

\section{Passage with $\triangle t\rightarrow 0^{+}$}

In this section we wish to present the proof of Theorem \ref{theorem2}, i.e. to demonstrate the passage with $\triangle t\rightarrow 0^{+}$ . The two previous section enable us to restrict attention to the case when the weak solution of the system (\ref{eq:system})-(\ref{slip}) exists, as provided by Theorem \ref{theorem1}. Our approach will be based on  some estimates uniform with respect to the length of time interval $\triangle t$ that we are going to gain here too. The task requires to work in the Bochner Spaces, but first let us introduce suitable notation:    
\begin{equation}\label{notacja}
\left.
\begin{array}{l}
\hat{\phi}(x,t)=\phi^{k}(x)\\
\tilde{\phi}(x,t)=\phi^{k}(x)+(t-k\triangle t)(\frac{\phi^{k+1}-\phi^{k}}{\triangle t})(x)
\end{array}
\right\}  \quad if \ k\triangle t\leq t<(k+1)\triangle t.
\end{equation}
This converts  our original system into 
\begin{equation}\label{sysweak}
\begin{array}{rl}
{}&\frac{\partial\tilde{\varrho}}{\partial t}+\dyv(\hat{\varrho}\hat{v})=0\quad in\ \Omega,\\
{}&\frac{\partial\widetilde{\varrho v}}{\partial t}+\dyv(\hat{\varrho}\hat{v}\varotimes \hat{v})-\mu\Delta \hat{v}-(\mu+\nu)\nabla\dyv \hat{v}+\nabla\pi(\hat{\varrho})=0\quad in\ \Omega,\\
{}&\hat{v}\cdot n=0\quad at\ \partial\Omega,\\
{}&n\cdotp T(\hat{v},\pi)\cdotp\tau+f\hat{v}\cdotp\tau=0\quad at\ \partial\Omega
\end{array}
\end{equation}
Moreover, the estimates (\ref{osz1}) and (\ref{vH1}) from the previous section now read:
\begin{eqnarray}
&\bullet&\hat{\varrho},\tilde{\varrho}\ are\  bounded \ in \  L_{\infty}(0,T;L_{\gamma}(\Omega))\label{i}\\
&\bullet&\hat{\varrho}\hat{v}^{2}, \widetilde{\varrho v^{2}}\ are\ bounded\ in\ L_{\infty}(0,T;L_{1}(\Omega))\label{j}\\
&\bullet&\hat{v},\ \tilde{v}\ are\ bounded\ in\ L_{2}(0,T;H^{1}(\Omega))\label{m}\\
&\bullet&\hat{\varrho}\hat{v},\widetilde{\varrho v}\ are\ bounded\ in\ L_{\infty}(0,T;L_{\frac{2\gamma}{\gamma+1}}(\Omega))\cup L_{2}(0,T;L_{r}(\Omega))\label{k}
\end{eqnarray} 
for $ 1\leq r<\gamma$ the last one holds as
\begin{equation*}
\|\varrho^{k}v^{k}\|_{2\gamma/(\gamma+1)}\leq\|\varrho^{k}\|^{1/2}_{\gamma}\|\varrho^{k}(v^{k})^{2}\|_{1}^{1/2}\quad and\quad\|\varrho^{k}v^{k}\|_{r}\leq\|\varrho^{k}\|_{\gamma}\|v^{k}\|_{\infty-\epsilon},
\end{equation*}
and all the bounds are independent of $\triangle t$. \\
Our next aim will be to reconstruct the estimation for the norm of pressure $\pi(\hat{\varrho})=\hat{\varrho}^{\gamma}$ in $L_{q}(\Omega\times(0,T))$ for some $q>1$. Unfortunately, as we have  seen in (\ref{2gamma}), such an estimate might not be true while $q=2$, but it turns out to work for $q=1+(1/\gamma)$. To show this we test the momentum equation with a function $\Phi$ of the form:
\begin{eqnarray*}
\Phi^{k}=\mathcal{B}((\varrho^{k})-\{\varrho^{k}\}),\quad in\ \Omega\\
\Phi^{k}=0\quad at\ \partial\Omega
\end{eqnarray*}
From this testing we obtain the following identity:
\begin{multline*}
\calka(\varrho^{k})^{\gamma+1}=\calka(\varrho^{k})^{\gamma}\{\varrho^{k}\}-\calka
\varrho^{k}v^{k}\varotimes v^{k}:\nabla\Phi^{k}+\mu\calka\nabla v^{k}:\nabla\Phi^{k}+(\mu+\nu)\calka\dyv v^{k}\dyv\Phi^{k}\\
+\calka\dt(\varrho^{k}v^{k}-\varrho^{k-1})\Phi^{k}=\sum_{i=1}^{5}I_{i}.
\end{multline*}
Multiplying by $\triangle t$, summing over $k=1,\ldots, M$ and employing our notation we get
\begin{multline}
\calkat \calka \hat{\varrho}^{\gamma+1}=\calkat \calka(\hat{\varrho})^{\gamma}\{\hat{\varrho}\}-\calkat\calka
\hat{\varrho}\hat{v}\varotimes \hat{v}:\nabla\hat{\Phi}+\mu\calkat\calka\nabla \hat{v}:\nabla\hat{\Phi}+(\mu+\nu)\calkat\calka\dyv \hat{v}\dyv\hat{\Phi}\\
+\calkat\calka\dt(\hat{\varrho}\hat{v}-\hat{\varrho}(\cdot-\triangle t)\hat{v}(\cdot-\triangle t))\hat{\Phi}=\sum_{i=1}^{5}I_{i}.
\end{multline}
We go one with estimations for each of terms separately.\\
(i) Since $\hat{\varrho}$ is bounded in $L_{\infty}(L_{1})$ and $L_{\infty}(L_{\gamma})$ one gets
\begin{equation*}
I_{1}=\calkat \calka(\hat{\varrho})^{\gamma}\{\hat{\varrho}\}=\calkat\frac{1}{|\Omega|}\|\hat{\varrho}\|_{L_{1}(\Omega)}\|\hat{\varrho}\|_{L_{\gamma}(\Omega)}^{\gamma}\leq CT.
\end{equation*}
(ii) The H$\mathrm{\ddot{o}}$lder's inequality, (\ref{m}) and (\ref{k}) imply
\begin{equation*}
I_{2}=-\calkat\calka
\hat{\varrho}\hat{v}\varotimes \hat{v}:\nabla\hat{\Phi}\leq\calkat\|\hat{v}\hat{\varrho}\|_{2\gamma/(\gamma+1)}\|\hat{v}\|_{W^{1}_{2}}\|\nabla\hat{\Phi}\|_{\gamma+1}\leq CT^{(\gamma-1)/2(\gamma+1)}\|\hat{\varrho}\|_{L_{\gamma+1}(L_{\gamma+1})}.
\end{equation*}
(iii) Due to the properties of tha Bogovskii functional $\|\nabla\Phi^{k}\|_{p}\leq c(p,\Omega)\|\varrho^{k}\|_{p}$, thus
\begin{multline*}
I_{3}+I_{4}=\mu\calkat\calka\nabla \hat{v}:\nabla\hat{\Phi}+(\mu+\nu)\calkat\calka\dyv \hat{v}\dyv\hat{\Phi}\leq\calkat\|\nabla\hat{v}\|_{L_{2}}\|\nabla\hat{\Phi}\|_{L_{\gamma+1}}\\
\leq CT^{(\gamma-1)/2(\gamma+1)}\|\hat{\varrho}\|_{L_{\gamma+1}(L_{\gamma+1})}.
\end{multline*}
(iv)
By the assumption that $\gamma>2$ we know that $\widehat{\widetilde{\varrho v}}\in L_{2}(0,T;L_{2}(\Omega))$ which is the special case of (\ref{k}), hence by the continuity equation
\begin{multline*}
I_{5}=\calkat\calka\dt(\hat{\varrho}\hat{v}-\hat{\varrho}(\cdot-\triangle t)\hat{v}(\cdot-\triangle t))\hat{\Phi}\\=\calkat\calka\frac{\partial}{\partial t}\widetilde{\varrho v \Phi}+\calkat\calka\dt\hat{\varrho}(\cdot-\triangle t)\hat{v}(\cdot-\triangle t)(\hat{\Phi}(\cdot-\triangle t)-\hat{\Phi})\\
\leq \sup_{0\leq t\leq T}\calka |\widetilde{\varrho v \Phi}|+\calkat \|\hat{\varrho}(\cdot-\triangle t)\hat{v}(\cdot-\triangle t)\|_{L_{2}(\Omega)}\|\hat{\varrho}(t)\hat{v}(t)\|_{L_{2}(\Omega)}\\
\leq C+\calkat \|\hat{\varrho}\|_{L_{\gamma}}^{2}\|\hat{v}\|^{2}_{L_{2\gamma/(\gamma-2)}}\leq C\end{multline*}
All together leads to desired conclusion
\begin{equation*}
\|\hat{\varrho}\|^{\gamma+1}_{L_{\gamma+1}(L_{\gamma+1})}\leq C\left( 1+T+T^{(\gamma-1)/2(\gamma+1)}\|\hat{\varrho}\|_{L_{\gamma+1}(L_{\gamma+1})} \right),
\end{equation*}
in particular we have
\begin{equation}
\sum_{k=1}^{M}\triangle t \|\varrho^{k}\|^{\gamma+1}_{L_{\gamma+1}}<C(T)\label{l}.
\end{equation}

We are now in a position to validate that  as $\triangle t\rightarrow 0$
the following convergences hold:
\begin{equation}\label{zbg}
[\hat{\varrho}-\hat{\varrho}(\cdot-\triangle t)],[\hat{\varrho}-\tilde{\varrho}]\rightarrow 0\quad
 in \ L_{q}(L_{\gamma})
 \end{equation}
for $q\in[1,\infty)$
\begin{equation}\label{zbg1}
[\hat{\varrho}\hat{v}-\hat{\varrho}\hat{v}(\cdot-\triangle t)],\ [\hat{\varrho}\hat{v}-\widetilde{\varrho v}]\rightarrow 0\quad
 in \ L_{q}(L_{r}),
 \end{equation}
for $\{q\in[1,\infty),\ r\in[1,\frac{2\gamma}{\gamma+1}]\}\cup\{q\in[1,2),\ r\in[1,\gamma)\}$,
\begin{equation}\label{zbg2}
[\hat{\varrho}\hat{v}\varotimes\hat{v}-\widetilde{\varrho v} \varotimes\hat{v}]\rightarrow 0\quad in \ L_{1}(L_{r})\cap L_{q}(L_{1}),
\end{equation}
for $q\in[1,\infty)\  r\in[1,\gamma)$.\\
To see this it suffices to use the estimates (\ref{i}, \ref{j}, \ref{m}, \ref{k}) together with the observation derived from (\ref{osz3}), namely 
\begin{equation}\label{zb_varrho}
\|\hat{\varrho}-\hat{\varrho}(\cdot-\triangle t)\|^{\gamma}_{L_{\gamma}(L_{\gamma})}\leq\triangle t C,
\end{equation}
moreover for the remaining term in (\ref{osz2}) we also have
\begin{equation}
\|\hat{\varrho}|\hat{v}-\hat{v}(\cdot-\triangle t)|^{2}\|_{L_{1}(L_{1})}\leq\triangle t C.
\end{equation}
From what has already been written we deduce that 
\begin{eqnarray}
\hat{\varrho},\ \tilde{\varrho}&\rightharpoonup&\varrho\quad weakly^{*}\ in\ L_{\infty}(L_{\gamma}),\ weakly\  in\  L_{\gamma+1}((0,T)\times\Omega),\label{uno}\\
\hat{v}&\rightharpoonup& v\quad weakly\ in\ L_{2}(H^{1}).\label{dos}
\end{eqnarray}
\begin{remark}\label{rem1}
Since  $\tilde{\varrho}\ \hat{\varrho},\ \hat{v}$ satisfy continuity equation $(\ref{sysweak})_{1}$, thus the sequence of functions $f(t)=\left(\calka\tilde{\varrho}\phi\ dx \right)(t)$ is bounded and equicontinuous in $C[0,T]$ for all $\phi\in C^{\infty}(\overline{\Omega}),\ \phi\cdot n=0$ at $\partial\Omega$. Therefore, the Arzela-Ascoli theorem, the density argument and the convergence established in (\ref{zbg}) yield the following
\begin{equation}
\hat{\varrho},\  \tilde{\varrho}\rightharpoonup \varrho\quad
 in \ C_{weak}(L_{\gamma}).
\end{equation}
\end{remark}

\noindent What is left is to show that we also have the corresponding  convergence of the products $\hat{\varrho} \hat{v},\ \hat{\varrho} \hat{v}\varotimes \hat{v}$. This can be done by repeated application of the following lemma.
\begin{lemat}\label{LLions}
Let $g^{n},\ h^{n}$ converge weakly to $g,\ h$ respectively in $L_{p_{1}}(L_{p_{2}}), \ L_{q_{1}}(L_{q_{2}})$ where $1\leq p_{1},p_{2}\leq\infty$ and
\begin{equation*}
\frac{1}{p_{1}}+\frac{1}{q_{1}}=\frac{1}{p_{2}}+\frac{1}{q_{2}}=1.
\end{equation*}
Let assume in addition that
\begin{equation}\label{condI}
\frac{\partial g^{n}}{\partial t}\ is\ bounded\ in\ L_{1}(W^{-m}_{1})\ for\ some\ m\geq0\ independent\ of\ n
\end{equation}
\begin{equation}\label{condII}
\|h^{n}-h^{n}(\cdot+\xi,t)\|_{L_{q_{1}}(L_{q_{2}})}\rightarrow 0\ as\ |\xi|\rightarrow 0,\ uniformly\ in\ n.
\end{equation}
Then $g^{n}h^{n}$ converges to $gh$ in the sense of distributions on $\Omega\times(0,T)$.
\end{lemat}
\noindent For the proof we refer the reader to \cite{PLL}. \\
\\
For our case, since $\frac{\partial\widetilde{\varrho}}{\partial t}$ is bounded in $L_{\infty}(W^{-1}_{2\gamma/(\gamma+1)})$
and $\frac{\partial\widetilde{\varrho v}}{\partial t}$ is bounded in $L_{\infty}(W^{-1}_{1})+L_{2}(H^{-1})$, the condition (\ref{condI}) is satisfied for $g^{n}=\widetilde{\varrho}, \widetilde{\varrho v}$ and $m=1$ respectively. Additionally, we have that since $h^{n}=\hat{v}$ is bounded in $L_{2}(H^{1})$  the condition (\ref{condII}) also holds true. \\
With this manner we see that $\widetilde{\varrho}\hat{v}$ converges weakly/weakly$^{*}$ in $L_{\infty}(L_{2\gamma/(\gamma+1)})$ and in $ L_{2}(L_{r})$ for $r\in[1,\gamma)$ to $\varrho v$ and that $\widetilde{\varrho v}\varotimes\hat{v}$ converges weakly in $L_{1}(L_{r})\cap L_{q}(L_{1})$,
for $q\in[1,\infty)\  r\in[1,\gamma)$ to $\varrho v\varotimes v$. Thus, the relations  (\ref{zbg1}) and (\ref{zbg2}) cause that we actually have
\begin{equation}
\hat{\varrho} \hat{v} \rightharpoonup\varrho v\quad  weakly\  in\ L_{q}(L_{r})
\end{equation}
for $\{q\in[1,\infty),\ r\in[1,\frac{2\gamma}{\gamma+1}]\}\cup\{q\in[1,2),\ r\in[1,\gamma)\}$,
\begin{equation}
\hat{\varrho} \hat{v}\varotimes \hat{v}\rightharpoonup\varrho v\varotimes v\quad weakly\ in \ L_{1}(L_{r})\cap L_{q}(L_{1}),
\end{equation}
for $q\in[1,\infty)\  r\in[1,\gamma)$.\\
Having this we can pass to the (weak,weak*) limit as $\triangle t\rightarrow 0^{+}$ in the system (\ref{sysweak}) everywhere expect in the term corresponding to
the pressure:
\begin{equation}\label{systembar}
\begin{array}{rl}
{}&\frac{\partial\varrho}{\partial t}+\dyv(\varrho v)=0\quad in\ \Omega,\\
{}&\frac{\partial\varrho v}{\partial t}+\dyv(\varrho v\varotimes v)-\mu\Delta v-(\mu+\nu)\nabla\dyv v+\nabla\overline{\pi(\varrho)}=0\quad in\ \Omega,\\
{}&v\cdot n=0\quad at\ \partial\Omega,\\
{}&n\cdotp T(v,\pi)\cdotp\tau+fv\cdotp\tau=0\quad at\ \partial\Omega
\end{array}
\end{equation}
From now on we will be using the following denotation
\begin{eqnarray*}
\mathbb{S}(\nabla v)&=&\mu(\nabla v+\nabla^{\bot}v)+\nu\dyv_{x}vI,\\
\mathbb{S}(\nabla \hat{v})&=&\mu(\nabla \hat{v}+\nabla^{\bot}\hat{v})+\nu\dyv_{x}\hat{v}I.
\end{eqnarray*}
The proof of strong convergence of $\pi(\varrho^{k})=(\varrho^{k})^{\gamma}$ in $L_{1}(\Omega\times(0,T))$ is based on some properties of the double Riesz transform, defined on the whole $\mathbb{R}^{2}$ in the following way
\begin{equation*}\mathcal{R}_{i,j}=-\partial_{x_{i}}(-\Delta)^{-1}_{x}\partial_{x_{j}},
\end{equation*}
where the inverse Laplacian is identified through the Fourier transform $\mathcal{F}$ and the inverse Fourier transform $\mathcal{F}^{-1}$ as
\begin{equation*}
(-\Delta)^{-1}(v)=\mathcal{F}^{-1}\left(\frac{1}{|\xi|^{2}}\mathcal{F}(v)\right).
\end{equation*}
We will be using general results on such operators as continuity but also some facts concerning the commutators involving Riesz operators, being mostly a consequence of Div-Curl lemma \cite{TT} or that of Coifman-Mayer \cite{CM}, \cite{FNP}. The best overall reference here for both: auxiliary tools and the general idea of the proof is \cite{FN}.\\
To take advantage of what we mentioned, there is a need to extended the system (\ref{sysweak}) to the whole $\mathbb{R}^{2}$, as this is where the definition of the operator $\Delta_{x}^{-1}$ makes sense.
We first observe that it can easily be done so for the continuity equation as $\hat{\varrho}\hat{v}\cdot n=0$ at $\partial\Omega$, hence
\begin{equation}
\frac{\partial1_{\Omega}\tilde{\varrho}}{\partial t}+\dyv(1_{\Omega}\hat{\varrho}\hat{v})=0.
\end{equation}
For the momentum equation  $(\ref{sysweak})_{2}$ we check that
\begin{eqnarray*}
\hat{\varphi}(t,x)=\psi(t)\zeta(x)\tilde{\phi},\quad \tilde{\phi}=(\nabla_{x}\Delta_{x}^{-1})[1_{\Omega}\tilde{\varrho}],\\ \psi\in C^{\infty}_{c}((0,T)),\ \zeta\in C^{\infty}_{0}(\overline{\Omega}),
\end{eqnarray*}
is an admissible test function. This can be seen as a consequence of  estimates (\ref{i}, \ref{j}, \ref{m}, \ref{k}, \ref{l}) and by the fact that the operator $\nabla_{x}\Delta_{x}^{-1}$ gives rise to the spatial regularity to its range comparing to its argument of one. Particularly, later on we will take advantage of that for $\gamma>2$, the embedding $W^{1}_{\gamma}(\Omega)\subset C(\overline{\Omega})$ together with Remark \ref{rem1} imply
\begin{equation}\label{CCC}
(\nabla_{x}\Delta_{x}^{-1})[1_{\Omega}\tilde{\varrho}]\rightarrow(\nabla_{x}\Delta_{x}^{-1})[1_{\Omega}\varrho]\quad in\ C([0,T]\times \overline{\Omega}).
\end{equation}
Having disposed of this preliminary step, we can get the following integral identity 
\begin{equation}\label{przed}
\calkat\calka\psi\zeta\left(\hat{\varrho}^{\gamma}\tilde{\varrho}+\mathbb{S}(\nabla\hat{v}):\nabla_{x}\Delta_{x}^{-1}\nabla_{x}[1_{\Omega}\tilde{\varrho}]\right)dx\ dt=\sum_{i=1}^{5}I_{i}
\end{equation}
where
\begin{eqnarray*}
I_{1}&=&\calkat\calka\psi\zeta\left(\widetilde{\varrho v}\Dt\tilde{\phi} +\hat{\varrho}\hat{v}\varotimes\hat{v}:\nabla_{x}\Delta_{x}^{-1}\nabla_{x}[1_{\Omega}\tilde{\varrho}]\right)\ dx\ dt,\\
I_{2}&=&-\calkat\calka\psi\hat{\varrho}^{\gamma}\nabla_{x}\zeta\cdot\nabla_{x}\Delta_{x}^{-1}[1_{\Omega}\tilde{\varrho}]\ dx \ dt,\\
I_{3}&=&\calkat\calka\psi \mathbb{S}(\nabla\hat{v}):\nabla_{x}\zeta\varotimes\nabla_{x}\Delta_{x}^{-1}[1_{\Omega}\tilde{\varrho}]\ dx\ dt,\\
I_{4}&=&-\calkat\calka\psi\left(\hat{\varrho}\hat{v}\varotimes\hat{v}\right):\nabla_{x}\zeta\varotimes\nabla_{x}\Delta_{x}^{-1}[1_{\Omega}\tilde{\varrho}]\ dx\ dt,\\
I_{5}&=&-\calkat\calka\Dt\psi\zeta\widetilde{\varrho v}\cdot\nabla_{x}\Delta_{x}^{-1}[1_{\Omega}\tilde{\varrho}]\ dx\ dt.
\end{eqnarray*}
Analogically, if we test the limit momentum equation by the corresponding test function
\begin{equation}
\varphi(t,x)=\psi(t)\zeta(x)\phi,\quad \phi=(\nabla_{x}\Delta_{x}^{-1})[1_{\Omega}\tilde{\varrho}],\ \psi\in C^{\infty}_{c}((0,T)),\ \zeta\in C^{\infty}_{0}(\overline{\Omega}),
\end{equation}
we get
\begin{equation}\label{po}
\calkat\calka\psi\zeta\left(\overline{\varrho^{\gamma}}\varrho+\mathbb{S}(\nabla v):\nabla_{x}\Delta_{x}^{-1}\nabla_{x}[1_{\Omega}\varrho]\right)dx\ dt=\sum_{i=1}^{5}I_{i}
\end{equation}
where
\begin{eqnarray*}
I_{1}&=&\calkat\calka\psi\zeta\left(\varrho v\Dt\phi +\varrho v\varotimes v:\nabla_{x}\Delta_{x}^{-1}\nabla_{x}[1_{\Omega}\varrho]\right)\ dx\ dt,\\
I_{2}&=&-\calkat\calka\psi\overline{\varrho^{\gamma}}\nabla_{x}\zeta\cdot\nabla_{x}\Delta_{x}^{-1}[1_{\Omega}\varrho]\ dx \ dt,\\
I_{3}&=&\calkat\calka\psi \mathbb{S}(\nabla v):\nabla_{x}\zeta\varotimes\nabla_{x}\Delta_{x}^{-1}[1_{\Omega}\varrho]\ dx\ dt,\\
I_{4}&=&-\calkat\calka\psi\left(\varrho v\varotimes v\right):\nabla_{x}\zeta\varotimes\nabla_{x}\Delta_{x}^{-1}[1_{\Omega}\varrho]\ dx\ dt,\\
I_{5}&=&-\calkat\calka\Dt\psi\zeta\varrho v\cdot\nabla_{x}\Delta_{x}^{-1}[1_{\Omega}\varrho]\ dx\ dt.
\end{eqnarray*}
The observation (\ref{CCC}) together with the consequences of lemma \ref{LLions} justify the convergences of the integrals $I_{2},\ldots, I_{5}$ from (\ref{przed}) to their counterparts in (\ref{po}). Thus we are left with the following identity
\begin{multline}\label{karaluch}
\lim_{\Delta t\rightarrow 0}\calkat\calka\psi\zeta\left(\hat{\varrho}^{\gamma}\tilde{\varrho}-\mathbb{S}(\nabla\hat{v}):\mathcal{R}[1_{\Omega}\tilde{\varrho}]\right)\\-\lim_{\Delta t\rightarrow 0}\calkat\calka\psi\zeta\left(\widetilde{\varrho v}\Dt\tilde{\phi} -\hat{\varrho}\hat{v}\varotimes\hat{v}:\mathcal{R}[1_{\Omega}\tilde{\varrho}]\right)\ dxdt\\
=\calkat\calka\psi\zeta\left(\overline{\varrho^{\gamma}}\varrho-\mathbb{S}(\nabla v):\mathcal{R}[1_{\Omega}\varrho]\right)dxdt\\-\calkat\calka\psi\zeta\left(\varrho v\Dt\phi -\varrho v\varotimes v:\mathcal{R}[1_{\Omega}\varrho]\right)\ dxdt.
\end{multline}
By the continuity equation we obtain 
\begin{equation*}
\Dt\phi=\mathcal{R}[1_{\Omega}\varrho  v],
\end{equation*}
and the same we have for the test function in the approximate case, thus (\ref{karaluch}) may be rewritten as
\begin{multline}\label{karaluch1}
\lim_{\Delta t\rightarrow 0}\calkat\calka\psi\zeta\left(\hat{\varrho}^{\gamma}\tilde{\varrho}-\mathbb{S}(\nabla\hat{v}):\mathcal{R}[1_{\Omega}\tilde{\varrho}]\right)dxdt\\-\lim_{\Delta t\rightarrow 0}\calkat\calka\psi\zeta\left(\widetilde{\varrho v}\mathcal{R}[1_{\Omega}\hat{\varrho}  \hat{v}] -\hat{\varrho}\hat{v}\varotimes\hat{v}:\mathcal{R}[1_{\Omega}\tilde{\varrho}]\right)\ dxdt\\
=\calkat\calka\psi\zeta\left(\overline{\varrho^{\gamma}}\varrho-\mathbb{S}(\nabla v):\mathcal{R}[1_{\Omega}\varrho]\right)\ dxdt\\-\calkat\calka\psi\zeta\left(\varrho v\mathcal{R}[1_{\Omega}\varrho  v] -\varrho v\varotimes v:\mathcal{R}[1_{\Omega}\varrho]\right)\ dxdt.
\end{multline}
Now we will show that we actually have that the second terms on each of sides are equivalent as $\Delta t$ goes to $0$. For this purpose we will use the Feireisl's lemma \cite{FN} which is a consequence of the div-curl one.
\begin{lemat}
Let
\begin{eqnarray*}
V_{n}&\rightharpoonup& V\quad weakly\ in\ L_{p}(\mathbb{R}^{2}),\\
r_{n}&\rightharpoonup& r\quad weakly\ in\ L_{q}(\mathbb{R}^{2}),
\end{eqnarray*}
where
\begin{equation*}
\frac{1}{p}+\frac{1}{q}=\frac{1}{s}<1.
\end{equation*}
Then
\begin{equation*}
V_{n}\mathcal{R}(r_{n})-r_{n}\mathcal{R}(V_{2})\rightharpoonup V\mathcal{R}(r)-r\mathcal{R}(V)\quad weakly\ in\ L_{s}(\mathbb{R}^{2}).
\end{equation*}
\end{lemat}
We will apply this lemma to $r_{n}=\hat{\varrho}(t,\cdot),\ V_{n}=\hat{\varrho}\hat{v}(t,\cdot)$ after extending them by 0 on the rest of $\mathbb{R}^{2}$ and noticing that they satisfy assumptions of the lemma for $ p= 2\gamma/(\gamma+1),\ q=\gamma$. Therefore we can take $s=\frac{2\gamma}{3+\gamma}$ and thus, for a.a $t\in[0,T)$
\begin{equation*}
\hat{\varrho}\hat{v}\mathcal{R}(\hat{\varrho})(t)-\hat{\varrho}\mathcal{R}(\hat{\varrho}\hat{v})(t)\rightharpoonup\varrho v\mathcal{R}(\varrho)(t)-\varrho\mathcal{R}(\varrho v)(t)\quad weakly\  in\  L_{s}(\Omega)
\end{equation*}
if we aditionally assume that $\gamma>3$.\\
In view of this, the embedding $L_{\frac{2\gamma}{3+\gamma}}(\Omega)\subset W^{-1}_{2}(\Omega)$ and (\ref{dos}) we get that
\begin{multline*}
\lim_{\Delta t\rightarrow 0}\calkat\calka\psi\zeta\hat{v}\left(\hat{\varrho}\mathcal{R}[1_{\Omega}\hat{\varrho}  \hat{v}] -\hat{\varrho}\hat{v}\mathcal{R}[1_{\Omega}\hat{\varrho}]\right)\ dxdt\\=
\calkat\calka\psi\zeta v\left(\varrho \mathcal{R}[1_{\Omega}\varrho  v] -\varrho v\mathcal{R}[1_{\Omega}\varrho]\right)\ dxdt,
\end{multline*}
and the relations (\ref{zbg1}) and (\ref{uno}) allow us to reduce (\ref{karaluch1}) to
\begin{multline}\label{karaluch2}
\lim_{\triangle t\rightarrow 0}\calkat\calka\psi\zeta\left(\hat{\varrho}^{\gamma}\hat{\varrho}-\mathbb{S}(\nabla\hat{v}):\mathcal{R}[1_{\Omega}\hat{\varrho}]\right)dxdt\\
=\calkat\calka\psi\zeta\left(\overline{\varrho^{\gamma}}\varrho-\mathbb{S}(\nabla v):\mathcal{R}[1_{\Omega}\varrho]\right)\ dxdt.
\end{multline}
Now observe that by the fact that $\zeta\in  C^{\infty}_{0}(\overline{\Omega})$ we may integrate by parts the second term on the left hand side and we will get
\begin{multline}\label{karaluch2,5}
\calkat\calka\psi\zeta\ \mathbb{S}(\nabla\hat{v}):\mathcal{R}[1_{\Omega}\hat{\varrho}]\ dxdt= \calkat\calka\psi\mathcal{R}:\left[\zeta\ \mathbb{S}(\nabla\hat{v})\right]\hat{\varrho}\ dxdt\\
=\calkat\calka\psi(2\mu+\nu)\dyv \hat{v}\hat{\varrho}\ dxdt+\calkat\calka\psi\Big(\mathcal{R}:\left[\zeta\ \mathbb{S}(\nabla\hat{v})\right]\ -\ \zeta\ \mathcal{R}:\left[\mathbb{S}(\nabla\hat{v})\right]\Big)\hat{\varrho}\ dxdt.
\end{multline}
With the same manner we can transform the corresponding term in the limit on the right hand side of (\ref{karaluch2}). After passing with $\triangle t$ to the limit in (\ref{karaluch2,5}) we get
\begin{multline}
\lim_{\triangle t\rightarrow 0}\calkat\calka\psi\zeta\ \mathbb{S}(\nabla\hat{v}):\mathcal{R}[1_{\Omega}\hat{\varrho}]\ dxdt\\
=\calkat\calka\psi(2\mu+\nu)\overline{\dyv v \varrho}\ dxdt+\calkat\calka\psi\Big(\mathcal{R}:\left[\zeta\ \mathbb{S}(\nabla v)\right]\ -\ \zeta\ \mathcal{R}:\left[\mathbb{S}(\nabla v)\right],\Big) \varrho\ dxdt
\end{multline}
where the precise form of last term on the right is a consequence of Div-Curl lemma. Therefore (\ref{karaluch2}) reduces to
\begin{equation*}
\calkat\calka\psi\zeta\left(\overline{\varrho^{\gamma}\varrho}-\overline{\varrho\dyv_{x} v}\right)dxdt
=\calkat\calka\psi\zeta\left(\overline{\varrho^{\gamma}}\varrho-\varrho\dyv_{x} v\right)\ dxdt.
\end{equation*}
and since the choice of functions $\psi$ and $\zeta$ was arbitrary we have that:
\begin{equation}\label{karaluch3}
\overline{\varrho^{\gamma}\varrho}-\overline{\varrho\dyv_{x} v}
=\overline{\varrho^{\gamma}}\varrho-\varrho\dyv_{x} v.
\end{equation}
Next, we take $\delta>0$ and multiply the discrete version of the continuity equation by $\ln(\varrho^{k}+\delta)$. After integrating by parts over $\Omega$ one get
\begin{equation*}
\dt\calka(\varrho^{k}-\varrho^{k-1})\ln(\varrho^{k}+\delta)-\calka\varrho^{k}v^{k}\frac{\nabla\varrho^{k}}{\varrho^{k}+\delta}=0.
\end{equation*}
By the Lebesgue monotone convergence theorem we can pass with $\delta\rightarrow 0^{+}$ and then integrate by parts once more to find
\begin{equation*}
\dt\calka(\varrho^{k}-\varrho^{k-1})\ln(\varrho^{k})+\calka\dyv(v^{k})\varrho^{k}=0.
\end{equation*}
Recall that due to Theorem \ref{theorem1} we have  $\calka\varrho^{k}=\calka\varrho^{k-1}$, thus whereas $x\ln(x)$ is a convex function above equality may be changed into
\begin{equation}\label{potwora}
\dt\calka\left[\varrho^{k}\ln(\varrho^{k})-\varrho^{k-1}\ln(\varrho^{k-1})\right]dx+\calka\dyv(v^{k})\varrho^{k}\leq 0.
\end{equation}
Summing from $k=1$ to $k=M$, multiplying by $\Delta t$ and using the notation (\ref{notacja}) we transform  (\ref{potwora})  to
\begin{equation*}
\calka\widetilde{\varrho\ln(\varrho)}(T)\ dx
+\calkat\calka\hat{\varrho}\ \dyv_{x}\hat{v}\ dxdt\leq\calka\varrho\ln(\varrho)(0)\ dx,
\end{equation*}
thus after passing to the limit one get
\begin{equation}\label{biforek}
\calka\overline{\varrho\ln(\varrho)}(T)\ dx
+\calkat\calka\overline{\varrho\ \dyv_{x}v}\ dxdt\leq\calka\varrho\ln(\varrho)(0)\ dx,
\end{equation}
For the limit momentum equation, we take advantage of the fact that it is satisfied in the whole space in sense of distributions, thus the solution is automatically a renormalised solution, i.e. it is allowed to multiply the equation by $\ln(\varrho+\delta)$. Then we integrate over $\Omega$, pass to the limit as $\delta$ goes to $0^{+}$ and integrate with rescpect to time to get
\begin{equation}\label{afterek}
\calka\varrho\ln\varrho(T)\ dx+\calkat\calka\varrho\ \dyv_{x}v\ dx dt=\calka\varrho\ln\varrho(0)\ dx.
\end{equation}
By comparing the two results from (\ref{biforek}) and (\ref{afterek}) we get that
\begin{equation*}
\calka\overline{\varrho\ln(\varrho)}(T)\ dx
+\calkat\calka\overline{\varrho\ \dyv_{x}v}\ dxdt\leq\calka\varrho\ln\varrho(T)\ dx+\calkat\calka\varrho\ \dyv_{x}v\ dx dt.
\end{equation*}
As in the proof of previous theorem, by the 'standard arguments' we show that $\varrho\overline{\varrho^{\gamma}}\leq\overline{\varrho^{\gamma+1}}$ which together with (\ref{karaluch3}) provide
\begin{equation}
\calkat\calka\overline{\varrho\ \dyv_{x}v}\ dxdt\geq\calkat\calka\varrho\ \dyv_{x}v\ dx dt.
\end{equation}
The two last information joined give the desired information, namely
\begin{equation*}
\varrho\ln{\varrho}=\overline{\varrho\ln{\varrho}},
\end{equation*}
and finally, by the convexity of function $x\ln{x}$, we obtain
\begin{equation*}
\lim_{\triangle t\rightarrow 0^{+}}\hat{\varrho}=\varrho\ \quad a.e.\ in\ (0,T)\times\Omega\
\end{equation*}
that completes the proof of Theorem \ref{theorem2}.
\begin{flushright}
$\Box$
\end{flushright}
\begin{center}
ACKNOWLEDGEMENTS
\end{center}
The author wishes to express her gratitude to Piotr Bogusław Mucha and Milan Pokorn\'y for suggesting the problem, stimulating conversations and help during preparation of the paper.\\
The author was supported by the International Ph.D.
Projects Programme of Foundation for Polish Science operated within the
Innovative Economy Operational Programme 2007-2013 funded by European
Regional Development Fund (Ph.D. Programme: Mathematical
Methods in Natural Sciences) and partly supported by the MN Grant No N N201 547438.

\end{document}